\documentclass[11pt]{amsart}
\usepackage[utf8]{inputenc}
\usepackage{fixltx2e}
\usepackage{graphicx}
\usepackage{amscd}
\usepackage{amsmath}
\usepackage{amssymb, latexsym, amsthm, amsbsy, amsfonts, mathrsfs}
\usepackage{bm}

\usepackage{mathtools}

\newtheorem{thm}{Theorem}[section]

\newtheorem{lem}[thm]{Lemma}
\newtheorem{cor}[thm]{Corollary}

\newtheorem{defn}[thm]{Definition}

\newtheorem{fact}[thm]{Fact}

\newtheorem{exmp}[thm]{Example}

\newcommand{\al}{\alpha}
\newcommand{\be}{\beta}

\newcommand{\rig}{\rightarrow}
\newcommand{\mrig}{\mathrel{-\!\!\!\!\!\rightarrow}}

\newcommand{\Rig}{\Rightarrow}

\newcommand{\bcdw}{\mathbin{\boldsymbol\cdot}}
\newcommand{\bcdn}{\mbox{\boldmath{$\cdot$}}}

\newcommand{\seteq}{\mathrel{\mbox{\,:\!}=\nolinebreak }\,}

\newcommand{\sbA}{{\boldsymbol{A}}}

\newcommand{\sbC}{{\boldsymbol{C}}}
\newcommand{\sbD}{{\boldsymbol{D}}}

\newcommand{\sbL}{{\boldsymbol{L}}}

\newcommand{\sbS}{{\boldsymbol{S}}}

\newcommand{\sbX}{{\boldsymbol{X}}}

\newcommand{\Fma}{\boldsymbol{\mathit{Fm}}}
\newcommand{\Fms}{\mathit{Fm}}

\newcommand{\ov}{\overline}

\bmdefine{\Fm}{Fm}                                
\bmdefine{\A}{A}                                   
\bmdefine{\B}{B}
\bmdefine{\C}{C}
\bmdefine{\D}{D}

\bmdefine{\btau}{\tau}                                  
\bmdefine{\brho}{\rho}                                  

\newcommand{\diamondplus}{\rotatebox{45}{\scalebox{0.8}{$\boxtimes$}}}

\theoremstyle{definition}

\theoremstyle{remark}

\begin{document}
\title[The Algebraic Significance of Weak Excluded Middle Laws]{The Algebraic Significance of\\Weak Excluded Middle Laws}
\author{T.\ L\'{a}vi\v{c}ka}
\address{Institute of Information Theory and Automation of the Czech Academy of Sciences, Pod Vod\'{a}renskou v\v{e}\v{z}\'{i} 4, 182 00 Prague 8, Czech Republic}
\email{lavicka.thomas@gmail.com}
\author{T.\ Moraschini}
\address{Department of Philosophy, Faculty of Philosophy, University of Barcelona, Carrer de Montalegre 6, 08001, Barcelona, Spain}
\email{tommaso.moraschini@ub.edu}
\author{J.G.\ Raftery}
\address{Department of Mathematics and Applied Mathematics,
 University of Pretoria,
 Private Bag X20, Hatfield,
 Pretoria 0028,
South Africa}
\email{{james.raftery@up.ac.za}}
\keywords{Weak excluded middle law, inconsistency lemma, protoalgebraic logic, algebraizable logic, intuitionistic logic, modal logic, relevance logic\vspace{2mm}}
\thanks{The second author was supported by the research grant 2017 SGR 95 of the AGAUR from the Generalitat de Catalunya, by the I+D+i research project PID2019-110843GA-I00 \textit{La geometria de las logicas no-clasicas} funded by the Ministry of Science and Innovation of Spain, and by the Beatriz Galindo grant BEAGAL18/00040 funded by the Ministry of Science and Innovation of Spain.  The third author was supported in part by the National Research Foundation of South Africa (UID 85407).}


\makeatletter
\renewcommand{\labelenumi}{\text{(\theenumi)}}
\renewcommand{\theenumi}{\roman{enumi}}
\renewcommand{\theenumii}{\roman{enumii}}
\renewcommand{\labelenumii}{\text{(\theenumii)}}
\renewcommand{\p@enumii}{\theenumi(\theenumii)}
\makeatother

\begin{abstract}
For (finitary) deductive systems, we formulate a signature-inde\-pendent abstraction
of the \emph{weak excluded middle law} (WEML), which strengthens the existing general notion of an inconsistency lemma (IL).
Of special interest is the case where a quasivariety $\mathsf{K}$ algebraizes a deductive system $\,\vdash$.  We prove that, in this case, if $\,\vdash$ has a WEML (in the general sense) then
every relatively subdirectly irreducible member of $\mathsf{K}$ has a greatest proper $\mathsf{K}$--congruence;
the converse holds if $\,\vdash$ has an inconsistency lemma.
The result extends, in a suitable form, to all protoalgebraic logics.
A super-intuitionistic logic possesses a WEML iff it extends $\mathbf{KC}$.
We characterize the IL and the WEML for normal modal logics and for relevance logics.
A normal extension of $\mathbf{S4}$ has a global consequence relation with a WEML iff it extends $\mathbf{S4.2}$, while every axiomatic extension of $\mathbf{R^t}$ with an IL has a WEML.
\end{abstract}


\maketitle

\allowdisplaybreaks


\section{Introduction}

Jankov \cite{Jan68} proved in 1968 that the `weak excluded middle law' ${\neg p\vee\neg\neg p}$ axiomatizes the strongest super-intuitionistic logic
having the same
positive fragment
as the intuitionistic propositional calculus.
In the literature, this extension is called $\mathbf{KC}$, and it has several algebraic characterizations.
According to one of them, a variety $\mathsf{K}$ of Heyting algebras models
a logic that proves
${\neg p\vee\neg\neg p}$
iff every subdirectly irreducible member of $\mathsf{K}$
has a greatest proper congruence; cf.\ Gabbay \cite[Thm.\,19(a), p.\,67]{Gab81}.

The
main result
of the present paper
generalizes this
characterization of $\mathbf{KC}$ to a signature-independent framework.  It is in the spirit of the `bridge theorems' of abstract algebraic logic \cite{Cze01,Fon16} that correlate, for instance, syntactic interpolation or definability properties with model-theoretic amalgamation or epimorphism-surjectivity demands
{\cite{BH06,CP99,MRW20a}},
and deduction-like theorems with congruence extensibility properties
{\cite{BP88,BP,Cze01,Raf11}}.

Any such generalization must assume some properties of negation, in a suitably abstract form.
A familiar feature of intuitionistic and classical propositional logic is
that
\begin{equation}\label{il intuitionistic}
\textup{$\Gamma \cup \{\al\}$ is inconsistent iff \,$\Gamma\vdash\neg\al$.}
\end{equation}
In the classical case, there is a stronger variant,
which adds to (\ref{il intuitionistic})
that
\begin{equation}\label{il classical}
\textup{$\Gamma\cup\{\neg\al\}$ is inconsistent iff  \,$\Gamma\vdash\al$.}
\end{equation}

Signature-independent abstractions of (\ref{il intuitionistic}) and its conjunction with (\ref{il classical}) were formulated
in \cite{Raf13}, for finitary (but otherwise arbitrary) deductive systems $\,\vdash$, yielding general notions of an \emph{inconsistency lemma} and a \emph{classical inconsistency lemma}.
It emerged that, when some variety $\mathsf{K}$ algebraizes $\,\vdash$, then an
inconsistency lemma
amounts to the demand that the
finitely generated
congruences of members of $\mathsf{K}$ should form dually pseudo-complemented join semilattices, whereas a \emph{classical} inconsistency lemma signifies that $\mathsf{K}$ is filtral and its nontrivial members lack trivial subalgebras \cite{Raf13} (see \cite{CR17} for the case of quasivarieties).

Subsequently, L\'{a}vi\v{c}ka and P\v{r}enosil \cite{LP} observed that classical inconsistency lemmas can also be construed as abstract `excluded middle laws'.
Trading $\,\vdash\al\vee\neg\al$ in for the more flexible assertion
\begin{equation}\label{elm}
\textup{whenever \,$\Gamma\cup\{\al\}\vdash\beta$\, and \,$\Gamma\cup\{\neg\al\}\vdash\beta$, \,then \,$\Gamma\vdash\beta$,}
\end{equation}
their main notion conjoins abstractions of (\ref{il intuitionistic})
and (\ref{elm}).
It is proved in \cite{LP} that, when $\,\vdash$ is algebraized by a variety $\mathsf{K}$, then $\,\vdash$ has an excluded middle law of the abstract kind iff it has an (ordinary) inconsistency lemma and $\mathsf{K}$ is semisimple.

Of course, the analogue of (\ref{elm}) for $\,\vdash\neg\al\vee\neg\neg\al$ requires that
\begin{equation}\label{weml}
\textup{whenever \,$\Gamma\cup\{\neg\al\}\vdash\beta$\, and \,$\Gamma\cup\{\neg\neg\al\}\vdash\beta$, \,then \,$\Gamma\vdash\beta$.}
\end{equation}
Our general definition of a \emph{weak excluded middle law} (WEML) will be a signature-independent abstraction of the conjunction of (\ref{il intuitionistic}) and (\ref{weml}).

Suppose again (for simplicity) that $\,\vdash$ is algebraized by a variety $\mathsf{K}$.
%
In this setting, our main result states
that $\,\vdash$ has a WEML iff it has an inconsistency lemma and every subdirectly irreducible member of $\mathsf{K}$ has a greatest proper congruence (Theorem~\ref{rsi greatest}).  That characterization is invariant under category equivalence (Corollary~\ref{cat eq}), so the definition of a WEML is appropriately stable.

It follows from Theorem~\ref{rsi greatest} that a super-intuitionistic logic has a WEML (in the general sense) iff it proves ${\neg p\vee\neg\neg p}$.  When we restrict Theorem~\ref{rsi greatest} to the global consequence relations of normal extensions $\mathbf{L}$ of $\mathbf{S4}$, the `convergence axiom' $\Diamond\Box p\to\Box\Diamond p$ emerges as the counterpart of ${\neg p\vee\neg\neg p}$,
i.e., the systems of this kind with a WEML are just those for which $\mathbf{L}$
extends $\mathbf{S4.2}$ (the least modal companion of $\mathbf{KC}$).  In the context of relevance logics, we show that an axiomatic extension of $\mathbf{R^t}$ with an inconsistency lemma must have a WEML, and we characterize the extensions for which this is the case.

Actually, Theorem~\ref{rsi greatest} is formulated for quasivarieties, not only varieties, and we prove it in the still wider framework of protoalgebraic logics (Theorem~\ref{prsi greatest}).  The core of its proof is an argument concerning the structure of algebraic lattices and their semilattices of compact elements (Theorem~\ref{main theorem}).

\section{Preliminaries}\label{preliminaries section}

From now on, $\vdash$ denotes a fixed but arbitrary (sentential) {\em deductive system},
i.e., a substitution-invariant \emph{finitary} consequence
relation over formulas in some algebraic language, where the language comprises a signature and a fixed infinite set of variables.
(Finitarity is the demand that, whenever $\Gamma\vdash\al$, then $\Gamma'\vdash\al$ for some finite
$\Gamma'\subseteq\Gamma$.)
Among other standard abbreviations, we signify `$\Gamma\vdash\al$ for all $\al\in\Pi$'
by $\Gamma\vdash\Pi$, and `$\Gamma\vdash\Pi$ and $\Pi\vdash\Gamma$'
by $\Gamma\dashv\vdash\Pi$.

Algebras are assumed to have the type of $\,\vdash$, unless we say otherwise.
The universe of an algebra $\sbA$ is denoted as $A$, and is always assumed non-empty.

We assume a familiarity with the basic theory of deductive
systems and matrix semantics, cf.\ \cite{Cze01,Fon16,Woj88}.  If $\langle\sbA,F\rangle$
is a matrix model of $\,\vdash$, then $F$ is called a $\,\vdash$\,--{\em filter\/} of
the algebra $\sbA$.  Because the set of $\,\vdash$\,--filters of $\sbA$ is closed under arbitrary
intersections, it becomes a complete lattice when ordered by set inclusion.  This lattice
is algebraic (as $\,\vdash$ is finitary), so its compact elements are just the finitely generated
$\,\vdash$\,--filters of $\sbA$.  In $\sbA$, the
$\,\vdash$\,--filter
generated by a subset $Y$
is denoted as $\textup{Fg}^\sbA_{\,\vdash\,}Y$, while $F+^\sbA G$ stands for the
join of two $\,\vdash$\,--filters $F$ and $G$.

Recall that
$\,\vdash$\,--{\em theories\/} are just $\,\vdash$\,--filters of the absolutely free algebra $\Fma$
generated by the variables of $\,\vdash$, and
{\em substitutions\/} are
endomorphisms of $\Fma$.

Let $v_1,v_2,v_3,\dots$ be a
denumerable sequence of distinct variables of $\,\vdash$.  We sometimes abbreviate $v_1$ as $v$.
For each $n\in\mathbb{N}^+\seteq\{1,2,3,\dots\}$, we set
\[
\mathit{Fm}(n)=\{\be\in \mathit{Fm}: \text{\,the variables occurring in $\be$ are among $v_1,\dots,v_{n}$}\}.
\]
If $\xi\in \mathit{Fm}(n)$ and $\sbA$ is an algebra, with $a_1,\dots,a_{n}\in A$,
then
$\xi^\sbA(a_1,\dots,a_{n})$
denotes $h(\xi)$, where $h\colon\textbf{\em Fm}
\to\sbA$ is any homomorphism
such that $h(v_i)=a_i$ for $\mbox{$i=1,\dots, n$}$.
If $\Xi\subseteq \mathit{Fm}(n)$, then
\[
\text{$\Xi^\sbA(a_1,\dots,a_{n})$ \,abbreviates\,
$\{\xi^\sbA(a_1,\dots,a_{n}):\xi\in\Xi\}$.}
\]
We omit the superscripts in $\xi^\sbA$, $\Xi^\sbA$, $+^\sbA$ and $\textup{Fg}^\sbA_{\,\vdash}$ when $\sbA$ is
$\textbf{\em Fm}$.

\section{Inconsistency Lemmas}\label{il section}

A set $\Xi$ of formulas of $\,\vdash$ is said to be {\em
inconsistent in} $\,\vdash$ \,if $\Xi\vdash\al$ for all $\al\in \mathit{Fm}$.
%
%
%
%
%

Let $\Psi_n\subseteq \mathit{Fm}(n)$ for all $n\in\mathbb{N}^+$.  Following \cite{Raf13}, we call $\{\Psi_n:n\in\mathbb{N}^+\}$
an {\em IL-sequence\/} for $\,\vdash$ provided that,
whenever
$\Gamma\cup\{\al_1,\dots,\al_n\}\subseteq \mathit{Fm}$ (with $n\in\mathbb{N}^+$), then
\[
\textup{$\Gamma\cup\{\al_1,\dots,\al_n\}$ is inconsistent in
$\,\vdash$ \ iff \ $\Gamma\vdash\Psi_n(\al_1,\dots,\al_n)$.}
\]
In this case, for all $n\in\mathbb{N}^+$ and $\al_1,\dots,\al_n\in \mathit{Fm}$,
\begin{align*}\label{cil eq2}
& \Psi_n(\al_1,\dots,\al_n)\cup\{\al_1,\dots,\al_n\} \textup{ \,is inconsistent in $\,\vdash$,} \\
& \Psi_n(\al_1,\dots,\al_n)\dashv\vdash \Psi_n(\al_{f1},\dots,\al_{fn}) \textup{ \,for
any permutation $f$ of $1,\dots,n$,}
\end{align*}
and if
$\{\Phi_n:n\in\mathbb{N}^+\}$ is another IL-sequence for $\,\vdash$, then $\Psi_n\dashv\vdash\Phi_n$ for all $n$.\,\footnote{\,It can be shown that if $\,\vdash$ has an IL-sequence and $\Xi$ is inconsistent in $\,\vdash$, then so is $s[\Xi]$, for every substitution $s$,
but we shall not need to rely on this observation.}

An IL-sequence
$\{\Psi_n:n\in\mathbb{N}^+\}$ for $\,\vdash$ is said to be {\em elementary\/} if it consists of {\em finite\/} sets $\Psi_n$.

\begin{defn}\label{il def 2}
\textup{(\cite{Raf13})\,
We say that $\,\vdash$ has an {\em inconsistency lemma\/}---briefly an {\em IL\/}---if it has an
elementary IL-sequence.\,\footnote{\,This notion is referred to as a `finitary global IL' in \cite{LP} (where deductive systems are not assumed to be finitary); the two uses of the word `finitary' are unrelated.}
}
\end{defn}

When $\,\vdash$ has an IL-sequence $\{\Psi_n:n\in\mathbb{N}^+\}$, then it has an elementary IL-sequence $\{\Psi_n':n\in\mathbb{N}^+\}$ iff $\mathit{Fm}$ is compact in the lattice of $\,\vdash$\,--theories.  In this case, every algebra $\sbA$ has a greatest compact $\,\vdash$\,--filter, namely $A$, and
we can arrange that $\Psi_n'\subseteq\Psi_n$ for all $n\in\mathbb{N}^+$ (see \cite[Sec.\,3]{Raf13}).

\begin{exmp}\label{ipc example}
\textup{Intuitionistic and classical propositional logic have a common IL, which takes the form
\begin{equation*}
\textup{$\Gamma\cup\{\al_1,\dots,\al_n\}$ is inconsistent \ iff \
$\Gamma\vdash\neg(\al_1\wedge \,\dots\, \wedge\al_n)$,}
\end{equation*}
whereas
\[
\{\{v_1\rig (v_2\rig (\,\dots\,\rig (v_n\rig\bot)\dots))\}:n\in\mathbb{N}^+\}
\]
is an IL-sequence for
the $\rig,\bot$ fragment of intuitionistic logic.
The latter example illustrates the need to cater separately for different values of $n$ in the general definition of an IL-sequence.}
\end{exmp}

\begin{defn}\label{prot def}
\textup{(\cite{BP86,Cze85,Cze86})}\,
\textup{We say that $\,\vdash$ is
{\em protoalgebraic\/} if there exists $\mbox{$\Lambda\subseteq \mathit{Fm}(2)$}$
such that
$\,\vdash\Lambda(v_1,v_1)$ and
$\{v_1\}\cup \Lambda(v_1,v_2)\vdash v_2$.
(In this case, $\Lambda$ can be chosen finite, because $\,\vdash$ is finitary.)}
\end{defn}

\noindent
Numerous additional characterizations of protoalgebraicity are known, e.g., see \cite{Cze01,Fon16}.
If $\,\vdash$ is protoalgebraic and $v_1 \nvdash v_2$, then $\Lambda\neq\emptyset$, so no algebra has an empty $\,\vdash$\,--filter.
The process of filter generation in algebras is very complicated in general, but it improves
as follows in the protoalgebraic case:

\begin{lem}\label{prot lem}
\textup{(\cite[Prop.~6.12]{Fon16})}\,
Let\/ $\,\vdash$ be protoalgebraic, and let\/ $\sbA$ be an algebra, with\/ $Y\cup\{a\}\subseteq A$.

Then\/ $a\in \textup{Fg}^{\sbA\,}_{\,\vdash\,} Y$
iff there exist\/ $\Gamma\cup\{\al\}\subseteq \mathit{Fm}$
and a homomorphism\/ $h\colon\Fma
\to\sbA$ such that\/ $\Gamma\vdash\al$ and\/
$h[\Gamma]\subseteq Y\cup\textup{Fg}^{\sbA\,}_{\,\vdash\,}\emptyset$ and\/ $h(\al)=a$\textup{.}
\end{lem}

\begin{thm}\label{filgen thm3}
\textup{(\cite[Thm.~3.6]{Raf13})}\,
Let\/ $\{\Psi_{n}:n\in\mathbb{N}^+\}$ be an elementary IL-sequence for a proto\-algebraic deductive system\/
$\,\vdash$\textup{.}
Let\/ $F$ be a\/ $\,\vdash$\,--filter of an algebra\/ $\sbA$\textup{,} and let\/
$a_1,\dots,a_n\in A$\textup{,} where\/ $n\in\mathbb{N}^+$\textup{.}  Then
\begin{equation*}
\text{\quad$A=F+^\sbA\textup{Fg}^{\sbA\,}_{\,\vdash\,}\{a_1,\dots,a_n\}$
\,iff\/\, $\Psi_{n}^\sbA(a_1,\dots,a_n )\subseteq F
$\textup{.}}
\end{equation*}
\end{thm}

Here, Lemma~\ref{prot lem} is used in proving the forward implication.

\section{Dually Pseudo-Complemented Semilattices}

Let $\langle S;+\rangle$ be a
{\em join semilattice
with\/} $0$, i.e., an idempotent commutative semigroup
that
has a least element with respect to the order
\[
x\leq y\;\Longleftrightarrow\; x+y=y.
\]
For $a,b\in S$,
there is no guarantee that $a$ and $b$ have a greatest lower bound in $\langle S;\leq\rangle$, but we
abbreviate
\[
c\leq a \;\;\&\;\; c\leq b\;\;\&\;\;(\forall x\in S)((x\leq a \;\;\&\;\; x\leq b)\;\Longrightarrow\;x\leq c)
\]
as $a\bcdw b=c$, so that $\bcdn$ is a partial binary operation on $S$.

For $a,b\in S$, we call $b$ the \emph{dual pseudo-complement of\/} $a$ if $\langle S;\leq\rangle$ has a greatest element $1$, and
$b$ is the least element of $\langle S;\leq\rangle$ for which ${a+b=1}$.  In this case, we write $b=a^*$.

We say that $\langle S;+\rangle$ is {\em dually pseudo-complemented\/} if each of its elements has a dual pseudo-complement.  In this case,
$\langle S;\leq\rangle$ has a greatest element
and for all ${a,b\in S}$, we have $a^{**}\leq a$ and,
by \cite[(18)]{Fri62},
\begin{equation}\label{frink}
(a+b)^{**}=a^{**}+b^{**}.
\end{equation}

Recall that the compact elements of an algebraic lattice
always form a join semilattice with $0$, where $+$ is the inherited join operation.

\begin{lem}\label{1 lemma}
Let $\sbL=\langle L;\leq\rangle$ be an algebraic lattice whose join-semilattice\/ $\sbS=\langle S;+\rangle$
of compact elements is dually pseudo-complemented.  Let\/ $0$ and\/ $1$ be the least and greatest elements of $\sbL$, respectively.  Let $a\in L\backslash\{1\}$.  Then\/
\begin{enumerate}
\item\label{1 lemma one}
$1\in S$.

\item\label{1 lemma two}
If\/ $1$ is join-irreducible in the interval\/ \textup{$[a,1]\seteq\{d\in L:a\leq d\}$,} then the interval\/ \textup{$[a,1)\seteq\{d\in L:a\leq d<1\}$} has a greatest element.

\item\label{*}
Whenever $c+d=1,$ with $c\in S$ and $d\in L$, then $c^*\leq d$.
\end{enumerate}
Thus, for each $c\in S$, the dual pseudo-complement of $c$ in $\sbS$ is also the dual pseudo-complement of $c$ in $\sbL$.
\end{lem}
\begin{proof}
(\ref{1 lemma one})\,
As $0$ is compact in $\sbL$, the greatest element of $\sbS$ is $0^*$.  But $1$ is a join of elements of $\sbS$ (as $\sbL$ is algebraic), so $1=0^*\in S$.

(\ref{1 lemma two})\, The interval $[a,1]$
is a complete sublattice of $\sbL$.  Therefore, in this interval,
if $1$ is join-irreducible, then it
is completely join-irreducible
(because it is compact, by (\ref{1 lemma one})), whence $[a,1)$ has a greatest element.

(\ref{*})\,
Each $d\in L$ is a join of compact elements, so if $c+d=1$, with $c\in S$, then
$c+d'=1$ for some compact $d'\leq d$, as $1$ is compact.  Then $c^*\leq d'$, so $c^*\leq d$.
\end{proof}

The above discussion is relevant, because the $\,\vdash$\,--filter lattice of an algebra $\sbA$ is algebraic, and
the compact $\,\vdash$\,--filters of
$\sbA$ form a join semilattice with $0$ under the operation $+^\sbA$.
The semilattice order $\leq$ is just $\subseteq$, and $\textup{Fg}^\sbA_{\,\vdash}\,\emptyset$ is the least element.  The partial operation $\bcdn$ is therefore intersection (which need not be a total operation).  The connection between an IL and dual pseudo-complements, suggested by Theorem~\ref{filgen thm3}, is as follows.

\begin{thm}\label{dually pseudo-complemented thm2}
\textup{(\cite[Thm.~3.7]{Raf13})}\,
Let\/ $\,\vdash$ be a protoalgebraic
deductive system.
Then the following conditions are equivalent.
\begin{enumerate}
\item\label{pc1}
$\,\vdash$ has an inconsistency lemma.

\item\label{pc2}
For every algebra\/ $\sbA$\textup{,} the compact\/
$\,\vdash$\,--filters of\/ $\sbA$ form a dually pseudo-complemented
semilattice with respect to\/ $+^\sbA$\textup{.}

\item\label{pc3}
The join semilattice of compact\/ $\,\vdash$\,--theories is dually pseudo-comple\-mented.
\end{enumerate}
In this case, if\/ $\{\Psi_n:n\in\mathbb{N}^+\}$ is an elementary IL-sequence for\/ $\,\vdash$\textup{,} then for any $n\in\mathbb{N}^+$ and any elements $a_1,\dots,a_n$ of an algebra $\sbA$\textup{,} we have
\begin{align*}
& \left(\textup{Fg}^\sbA_{\,\vdash}\{a_1,\dots,a_n\}\right)^*=\textup{Fg}^\sbA_{\,\vdash}\,\Psi^\sbA_n(a_1,\dots,a_n).
\end{align*}
\end{thm}

The proof of Theorem~\ref{dually pseudo-complemented thm2} makes significant use of protoalgebraicity.

\section{Weak Excluded Middle Laws}

Suppose $\{\Psi_n:n\in\mathbb{N}^+\}$ is an elementary IL-sequence for $\,\vdash$.
For each $n\in\mathbb{N}^+$, let
$\#n=\left|\Psi_n\right|$ and $\Psi_n=\{\psi_n^1,\dots,\psi_n^{\#n}\}$, and define
\[
\Psi_{\#n}\Psi_n
\seteq
\Psi_{\#n}(\psi_n^1
,\dots,\psi_n^{\#n}
).
\]
(The \emph{definiens} is essentially unambiguous, by the remark on permutations in Section~\ref{il section}.)
Then, for any $\al_1,\dots,\al_n\in \mathit{Fm}$, we have
\begin{align*}\label{intuitionistic}
& \al_1,\dots,\al_n\vdash\Psi_{\#n}\Psi_n(\al_1,\dots,\al_n), \textup{ and}\\
& \Psi_{\#n}\Psi_n(\al_1,\dots,\al_n
)\cup \Psi_n(\al_1,\dots,\al_n
)
\textup{ is inconsistent in $\,\vdash$.}
\end{align*}
For any
elements $a_1,\dots,a_n$ ($n\in\mathbb{N}^+$) of an algebra $\sbA$\textup{,} Theorem~\ref{dually pseudo-complemented thm2} gives
\begin{align*}
& \left(\textup{Fg}^\sbA_{\,\vdash}\{a_1,\dots,a_n\}\right)^{**}=\textup{Fg}^\sbA_{\,\vdash}\,\Psi_{\#n}\Psi^\sbA_n(a_1,\dots,a_n)\textup{.}
\end{align*}

\begin{exmp}\label{ipc 2nd example}
\textup{The theorems of intuitionistic propositional logic ($\mathbf{IPC}$)
do not include the formula
${\neg v\vee\neg\neg v}$.
As we noted in the introduction, $\mathbf{KC}$ is the extension of $\mathbf{IPC}$ by the axiom $\neg v \vee \neg\neg v$.
When $\mathbf{L}$ is an axiomatic extension of $\mathbf{KC}$, the following implication holds (bearing the Deduction Theorem in mind):
\[
\textup{whenever \,$\Gamma,\neg\al\vdash_\mathbf{L}\be$\, and \,$\Gamma,\neg\neg\al\vdash_\mathbf{L}\beta$, \,then \,$\Gamma\vdash_\mathbf{L}\beta$.}
\]
This phenomenon is abstracted in the next definition.}
\end{exmp}
\begin{defn}\label{weml def}
\textup{We say that $\,\vdash$ has a \emph{weak excluded middle law} (WEML) if it has an elementary IL-sequence
$\{\Psi_n:n\in\mathbb{N}^+\}$ such that, for each $n\in\mathbb{N}^+$,
\[
\textup{if \,$\Gamma\cup\Psi_n(\al_1,\dots,\al_n)\vdash\varphi$ \,and\, $\Gamma\cup\Psi_{\#n}\Psi_n(\al_1,\dots,\al_n)\vdash \varphi$, then $\Gamma\vdash\varphi$.}
\]}
\end{defn}

We observed in Section~\ref{il section} that the form of an IL is unique up to inter-derivability.  It follows that if one elementary IL-sequence establishes a WEML for $\,\vdash$, then so does any other.  An IL persists in axiomatic extensions \cite[p.\,400]{Raf13}, and it is easy to see that the same applies to a WEML.  The presence of a WEML can be characterized as follows.

\begin{thm}\label{weak}
Let\/ $\,\vdash$ be a protoalgebraic deductive system.  Then the following conditions are equivalent.
\begin{enumerate}
\item\label{weak1}
$\,\vdash$ has a WEML.

\item\label{weak2}
For every algebra $\sbA$\textup{,} the join semilattice of\/ compact\/ $\,\vdash$\,--filters of\/ $\sbA$ is dually pseudo-complemented and satisfies\/
\[
(x+y^*)\bcdw (x + y^{**})=x.
\]
\item\label{weak3}
The join semilattice of\/ compact\/ $\,\vdash$\,--theories is dually pseudo-comple\-mented and satisfies\/
$(x+y^*)\bcdw (x + y^{**})=x$\textup{.}
\end{enumerate}
\end{thm}
\begin{proof}
In view of Theorem~\ref{dually pseudo-complemented thm2}, all three conditions imply that
$\,\vdash$ has an elementary IL-sequence $\{\Psi_n:n\in\mathbb{N}^+\}$, so let us assume this.
Let ${\Lambda\subseteq\mathit{Fm}(2)}$ be a finite set witnessing Definition~\ref{prot def}.
If $v_1\vdash v_2$, then (\ref{weak1})--(\ref{weak3}) are trivially true, so assume that $v_1 \nvdash v_2$.  Then $\Lambda\neq\emptyset$ and no algebra has an empty $\,\vdash$\,--filter.  Moreover, $\Psi_n\neq\emptyset$ for all $n\in\mathbb{N}^+$
(because $v_1\vdash\Psi_n(v_1,v_1,\dots,v_1)$ would entail $v_1\vdash v_2$, by the definition of an IL).

Trivially, (\ref{weak2}) implies (\ref{weak3}).
Theorem~\ref{dually pseudo-complemented thm2} shows that (\ref{weak3}) implies (\ref{weak1}),
remembering that $\,\vdash$ is finitary, that $\Gamma\vdash\beta$ paraphrases $\beta\in\textup{Fg}_{\,\vdash}\,\Gamma$, and that we always have $\textup{Fg}_{\,\vdash}(\Gamma\cup\Delta)=(\textup{Fg}_{\,\vdash}\,\Gamma)+(\textup{Fg}_{\,\vdash}\,\Delta)$.
It remains to prove that (\ref{weak1}) implies (\ref{weak2}).

Assuming (\ref{weak1}), let $G$ be a compact $\,\vdash$\,--filter of an algebra $\sbA$, and $\ov{b}\seteq b_1,\dots,b_n$ a finite sequence of elements of $A$, where $n\in\mathbb{N}^+$.
Let
$H=\textup{Fg}^\sbA_{\,\vdash}\,\{b_1,\dots,b_n\}$, so $H^*=\textup{Fg}^\sbA_{\,\vdash}\,\Psi_n^\sbA(\ov{b})$ and $H^{**}=\textup{Fg}^\sbA_{\,\vdash}\,\Psi_{\#n}\Psi_n^\sbA(\ov{b})$, by Theorem~\ref{dually pseudo-complemented thm2}.  Let $c\in (G+^\sbA H^*)\cap (G+^\sbA H^{**})$.  We need to show that $c\in G$.

By Lemma~\ref{prot lem} and the finitarity of $\,\vdash$, there exist a finite set
\[
\Pi=\Sigma_1\cup\Sigma_2\cup\Sigma_1'\cup\Sigma_2'\cup\{\varphi,\varphi'\}\subseteq \mathit{Fm}
\]
and homomorphisms $g,h\colon\Fma\to\sbA$ such that
\begin{center}
$\Sigma_1\cup\Sigma_2\vdash\varphi \textup{ \ and \ } \Sigma_1'\cup\Sigma_2'\vdash\varphi';$\\ \smallskip
$g[\Sigma_1]\cup h[\Sigma_1']\subseteq G \textup{ \ and \ } g[\Sigma_2]\subseteq \Psi_n^\sbA(\ov{b}) \textup{ \ and \ } h[\Sigma_2']\subseteq \Psi_{\#n}\Psi_n^\sbA(\ov{b});$\\ \smallskip
$g(\varphi)=c=h(\varphi').$
\end{center}
As $\Pi$ is finite,
the substitution-invariance of $\,\vdash$ allows us to assume, without loss of generality, that the variables occurring in members of $\Sigma_1\cup\Sigma_2\cup\{\varphi\}$ do not occur in any member of $\Sigma_1'\cup\Sigma_2'\cup\{\varphi'\}$.  We can therefore arrange that $g=h$, and also that $g(z_i)=b_i$ for $i=1,\dots,n$, where $\ov{z}=z_1,\dots,z_n$ is a sequence of distinct variables that are absent from all formulas in $\Pi$.

We construct a finite set $\Gamma\subseteq g^{-1}[G]$ as follows.  We stipulate that
\[
\Sigma_1\cup\Sigma_1'\subseteq\Gamma.
\]
Also, all elements of $\Lambda(\varphi,\varphi')\cup \Lambda(\varphi',\varphi)$ are included in $\Gamma$.  (These are sent by $g$ into $G$, because $\,\vdash\Lambda(v_1,v_1)$, and because $g(\varphi)=g(\varphi')$.)
Recall that $\Psi_n=\{\psi_n^1,\dots,\psi_n^{\#n}\}$, so
for each $\al\in\Sigma_2$, there exists $j$ such that $g(\al)=(\psi_n^j)^\sbA(\ov{b})$; we include all elements of $\Lambda(\psi_n^j(\ov{z}),\al)$ in $\Gamma$.
Likewise, as $\Psi_{\#n}\Psi_n=\Psi_{\#n}(\psi_n^1
,\dots,\psi_n^{\#n})$, we can choose,
for each $\al'\in\Sigma_2'$, a number $k$ so that $g(\al')=\psi_{\#n}^k(\psi_n^1,\dots,\psi_n^{\#n})^\sbA(\ov{b})$; we include all elements of
\[
\Lambda(\psi_{\#n}^k(\psi_n^1,\dots,\psi_n^{\#n})
(\ov{z}),\al')
\]
in $\Gamma$.  (All of these formulas belong to $g^{-1}[G]$, again since $\,\vdash\Lambda(v_1,v_1)$).  This completes the construction of $\Gamma$.

Because $\Sigma_1\cup\Sigma_2\vdash\varphi$, we have $\Gamma\cup\Psi_n(\ov{z})\vdash\varphi$ (by the rule
\[
\{v_1\}\cup\Lambda(v_1,v_2)\vdash v_2
\]
and the substitution-invariance and transitivity of $\,\vdash$).  Likewise, because $\Sigma_1'\cup\Sigma_2'\vdash\varphi'$, we have $\Gamma\cup\Psi_{\#n}\Psi_n(\ov{z})\vdash\varphi$.  Therefore, $\Gamma\vdash\varphi$, by (\ref{weak1}), and since $g[\Gamma]\subseteq G$, it follows that $c=g(\varphi)\in G$.
\end{proof}

Theorem~\ref{weak} implies
that, for a protoalgebraic deductive system with a WEML, the semilattice of compact deductive filters of any algebra must satisfy $y^*\bcdw y^{**}=0$.  A deductive system $\,\vdash$ is said to be \emph{filter-distributive} if every algebra has a distributive lattice of $\,\vdash$\,--filters.

\begin{cor}\label{weak cor}
Let\/ $\,\vdash$ be a filter-distributive
protoalgebraic
deductive system with an IL.  If the semilattice of compact $\,\vdash$\,--theories satisfies\/ $y^*\bcdw y^{**}=0$\textup{,} then\/ $\,\vdash$ has a WEML.
\end{cor}
\begin{proof}
Distributivity upgrades $y^*\bcdw y^{**}=0$ to $(x+y^*)\bcdw (x + y^{**})=x$.
\end{proof}
In fact, a protoalgebraic deductive system is filter-distributive iff it possesses a `parameterized disjunction' \cite[Thm.\,2.5.17]{Cze01} (also see \cite{CN13}).  As we shall not need to employ this syntactic notion, we omit its definition, which can be found in \cite[p.\,144]{Cze01}.

\section{Reduced Matrix Models}

Suppose $\Gamma\cup\{\varphi\}\subseteq\Fms$.  The following is well known (see \cite[Sec.\,3.7]{Woj88}).
\begin{fact}\label{completeness fact}
$\Gamma\vdash\varphi$ iff the implication
$(h[\Gamma]\subseteq F \:\Longrightarrow\: h(\varphi)\in F)$
holds for every homomorphism\/ $h\colon\Fma\to\sbA$ and every $\,\vdash$--filter $F$ of $\sbA$ such that
\begin{enumerate}
\item\label{cmi}
$F$ is completely meet-irreducible in the\/ $\,\vdash$\,--filter lattice of $\sbA,$ and
\item\label{reduced}
every congruence of $\sbA$ that identifies two distinct elements of $A$ also identifies an element of $F$ with a non-element of $F$.
\end{enumerate}
\end{fact}
\noindent
Condition~(\ref{reduced}) is more commonly phrased as `the matrix $\langle\sbA,F\rangle$ is \emph{reduced}'.  When (\ref{reduced}) is assumed, then (\ref{cmi}) may be rendered as `$\langle\sbA,F\rangle$ is $\,\vdash$\,--\emph{subdirectly irreducible}', because it means that $\langle\sbA,F\rangle$ cannot be decomposed subdirectly in the class of \emph{reduced} matrix models of $\,\vdash$.  Moreover, every reduced matrix model of $\,\vdash$ is isomorphic to a subdirect product of ones that are $\,\vdash$\,--subdirectly irreducible \cite[pp.\,242--3]{Woj88};
the finitarity of $\,\vdash$ is relied on here.

Fact~\ref{completeness fact} states that
the subdirectly irreducible reduced matrix models of $\,\vdash$ are always adequate as a semantics for $\,\vdash$.  In the absence of any data about $\,\vdash$, this is normally the most economical semantics at hand, and it takes the expected form in familiar examples.  (For classical propositional logic it yields just two-element Boolean algebras, with singleton filters comprising the top element in each case.)

We can now prove our main result about the WEML.  In fact, the proof can be carried out entirely in the setting of algebraic lattices.

\begin{thm}\label{main theorem}
Let $\sbL=\langle L;\leq\rangle$ be an algebraic lattice whose join-semilattice\/ $\sbS=\langle S;+\rangle$
of compact elements is dually pseudo-complemented.  Let\/ $0$ and\/ $1$ be the least and greatest elements of $\sbL$, respectively.  Then
the following conditions are equivalent\/\textup{:}
\begin{enumerate}
\item\label{(0)}
$a = (a + c^*) \bcdw (a + c^{**})$ for all $a,c \in S$\textup{;}

\item\label{(i)}
$a = (a + c^*) \bcdw (a + c^{**})$ for every $a \in L$ and $c \in S$\textup{;}

\item\label{(ii)}
whenever\/ $a \in L \backslash\{1\}$ is meet-irreducible in $\sbL,$ then the
interval\/ \textup{$[a, 1)$}
has a largest element\/\textup{;}

\item\label{(iii)}
whenever\/ $a \in L$ is completely meet-irreducible in $\sbL,$ then\/
$1$ is join-irreducible in the interval\/ \textup{$[a, 1]$.}
\end{enumerate}
\end{thm}
\begin{proof}
(\ref{(0)})$\:\Rig\:$(\ref{(i)}):
Let $a\in L$ and $c\in S$.  Let $b\in L$, with $b\leq a+c^*$ and $b\leq a+c^{**}$.  We need to show that $b\leq a$.  As $\sbL$ is algebraic, $b$ is the join (in $\sbL$) of some $\{b_i:i\in I\}\subseteq S$, and $a$ is also a join of compact elements of $\sbL$.  Let $i\in I$.  It suffices to show that $b_i\leq a$.  As $b_i\leq a+c^*,\, a+c^{**}$ and $b_i$ is compact, there is a compact element $a'$ of $\sbL$, with $a'\leq a$, such that $b_i\leq a'+c^*,\,a'+c^{**}$. So, because $a',c\in S$, it follows from (\ref{(0)}) that $b_i\leq a'$, whence $b_i\leq a$, as required.

(\ref{(i)})$\:\Rig\:$(\ref{(ii)}):
Let $a\in L\backslash\{1\}$ be meet-irreducible in $\sbL$.
We need to show that $1$ is completely join-irreducible
in the interval $[a, 1]$.
Suppose $1$ is the join of a family ${X \subseteq [a, 1]}$.  We must show that $1\in X$.
Since $\sbL$ is algebraic, each $x \in X$ is the join
of a family $Y_x$ of compact elements of $\sbL$.  Then $1$ is the join of \,$\bigcup_{x\in X}Y_x$, but
$1$ is itself compact (by Lemma~\ref{1 lemma}(\ref{1 lemma one})),
so $1$ is already the join of a finite subset $Y$ of \,$\bigcup_{x\in X}Y_x$.
Note that $Y\neq\emptyset$, because $1\neq 0$ (as
$a \in L \backslash \{1\}$).
Let $Y = \{y_1, \dots, y_n\}$, where $n$ is a positive integer, so $1=y_1+ {\,\dots\,} + y_n$.  As $Y$ consists of
compact elements, it follows from (\ref{(i)}) that
$a = (a + y_i^*) \bcdw (a + y_i^{**})$ for $i=1,\dots,n$.
Then, for each $i$, the meet-irreducibility of $a$ in $\sbL$ yields
$y_i^* \leq a$ or $y_i^{**}\leq a$.

If $y_i^*\leq a$ for some $i$, then $y_i+a=1$.  In this case, choose $x\in X$ with $y_i\in Y_x$, so $y_i,a\leq x$, whence $1=x\in X$, as required.  It therefore suffices to rule out the possibility that $y_i^*\nleq a$ for all $i$, i.e., that $y_i^{**}\leq a$ for all $i$.
Suppose, with a view to contradiction, that $y_1^{**}, \dots, y_n^{**} \leq a$.  Then
\[
1 = 1^{**} = (y_1 + {\,\dots\,} + y_n)^{**} = y_1^{**} + {\,\dots\,} + y_n^{**} \leq a,
\]
where the
third equality follows from (\ref{frink}).
But the displayed line contradicts the fact that $a \ne 1$, and this completes the proof of (\ref{(ii)}).\,\footnote{\,It would have sufficed here to show that $1$ is join-irreducible in $[a,1]$, owing to Lemma~\ref{1 lemma}(\ref{1 lemma two}), but that would not have simplified the argument significantly.}

(\ref{(ii)})$\:\Rig\:$(\ref{(iii)})  is clear, because $1$ is not completely meet-irreducible in $\sbL$ (as it is the meet of the empty subset of $L$).

(\ref{(iii)})$\:\Rig\:$(\ref{(0)}):
Let $a,c\in S$ and $b\in L$, with $b\leq a+c^*$ and $b\leq a+c^{**}$.  We need to show that $b\leq a$.
As $\sbL$ is algebraic, $a$ is the meet (in $\sbL$) of a set $\{a_i:i\in I\}$ of completely meet-irreducible elements of $\sbL$.
Let $i\in I$ be fixed.  It suffices to show that $b\leq a_i$.

We have $c+c^*=1$, hence $(a_i+c)+(a_i+c^*)=1$.  As $a_i$ is completely meet-irreducible,
$1$ is join-irreducible in $[a_i,1]$, by (\ref{(iii)}), whence $a_i+c$ or $a_i+c^*$ is $1$.

If $1=a_i+c$, then
$a,c^*\leq a_i$ (using Lemma~\ref{1 lemma}(\ref{*}) in the latter case), whence
$b\leq a+c^*\leq a_i$.

Likewise, if $1=a_i+c^*$, then $a,c^{**}\leq a_i$, whence $b\leq a_i$.
\end{proof}

Theorem~\ref{main theorem} specializes immediately as follows, in view of Theorem~\ref{weak} and Lemma~\ref{1 lemma}(\ref{1 lemma two}).

\begin{thm}\label{prsi greatest}
Let\/ $\,\vdash$ be a protoalgebraic deductive system with an inconsistency lemma.  Then the following conditions are equivalent.
\begin{enumerate}
\item\label{prsi greatest1}
$\,\vdash$ has a WEML.

\item\label{prsi greatest2}
Whenever\/ $F$ is a meet-irreducible\/ $\,\vdash$\,--filter of an algebra\/ $\sbA$\textup{,} with\/ $F\neq A$\textup{,} then the interval\/ $[F,A)$ in the\/ $\,\vdash$\,--filter lattice of $\sbA$ has a greatest element.

\item\label{prsi greatest3}
For each\/
$\,\vdash$\,--subdirectly irreducible reduced matrix model\/ $\langle\sbA,F\rangle$
of\/ $\,\vdash,$
the interval\/ $[F,A)$ of the\/ $\,\vdash$\,--filter lattice of $\sbA$
has a greatest element.
\end{enumerate}
\end{thm}

\section{Excluded Middle Laws}\label{excluded middle laws}

The following definition is due to L\'{a}vi\v{c}ka and P\v{r}enosil,
who formulated it in a wider framework.
\begin{defn}\label{eml def}
\textup{(\cite{LP})\, A deductive system $\,\vdash$ has an \emph{excluded middle law} (EML) if, for each $n\in\mathbb{N}^+$, there is a finite set $\Psi_n\subseteq \Fms(n)$ such that for all $\Gamma\cup\{\al_1,\dots,\al_n,\varphi\}\subseteq\Fms$,}
\begin{enumerate}
\item
\textup{$\Psi_n(\al_1,\dots,\al_n)\cup\{\al_1,\dots,\al_n\}$ is inconsistent in $\,\vdash$, and}\smallskip
\item
\textup{whenever $\Gamma\cup\{\al_1,\dots,\al_n\}\vdash\varphi$ and $\Gamma\cup\Psi_n(\al_1,\dots,\al_n)\vdash\varphi$, then $\Gamma\vdash\varphi$.}
\end{enumerate}
\end{defn}
\noindent
It is easy to see that, in this case, $\boldsymbol{\Psi}\seteq\{\Psi_n:n\in\mathbb{N}^+\}$ is an elementary IL-sequence for $\,\vdash$, and that $\boldsymbol{\Psi}$ also establishes a WEML for $\,\vdash$.  Theorem~\ref{weak} persists when we replace `WEML' by `EML', and $(x+y^*)\bcdw (x + y^{**})=x$ by
\[
(x+y)\bcdw (x + y^{*})=x
\]
in its statement; no essential change to its proof is needed.

For algebraic lattices, the EML-analogue of Theorem~\ref{main theorem} is as follows.

\begin{thm}\label{first cor}
Let $\sbL=\langle L;\leq\rangle$ be an algebraic lattice whose join-semilattice\/ $\sbS=\langle S;+\rangle$
of compact elements is dually pseudo-complemented.  Let\/ $0$ and\/ $1$ be the least and greatest elements of $\sbL$, respectively.  Then the following conditions are equivalent\/\textup{:}
\begin{enumerate}
\item\label{(i')}
$a = (a + c) \bcdw (a + c^{*})$ for every $a \in L$ and $c \in S$\textup{;}

\item\label{(ii')}
Whenever $a\in L\backslash\{1\}$ is meet-irreducible, then\/ $[a,1]=\{a,1\}$.

\item\label{(iii')}
Whenever $a\in L$ is completely meet-irreducible, then\/ $[a,1]=\{a,1\}$.
\end{enumerate}
\end{thm}
\begin{proof}
(\ref{(i')})$\:\Rig\:$(\ref{(ii')}):
Suppose $a<x\in L$.  We must show that $x=1$.  As $\sbL$ is algebraic, $x$ is the join of a set $W$ of elements of $S$.  As $x\nleq a$, there exists $c\in W$ such that $c\nleq a$, i.e., $a+c\neq a$.  As $c$ is compact, (\ref{(i')})
gives $a=(a+c)\bcdw(a+c^*)$, but $a$ is meet-irreducible, so $a=a+c^*$, i.e., $c^*\leq a$, whence $a+c=1$.  Now $a,c\leq x$, so $1=a+c\leq x$, i.e.,
$x=1$, as required.

(\ref{(ii')})$\:\Rig\:$(\ref{(iii')}) is clear.

(\ref{(iii')})$\:\Rig\:$(\ref{(i')}):  Let $a,b\in L$ and $c\in S$, where $b\leq a+c$ and $b\leq a+c^*$.  We must show that $b\leq a$.  As $\sbL$ is algebraic, $a$ is the meet of a family of completely meet-irreducible elements $a_i\in L$ ($i\in I)$.  Let $i\in I$.  It suffices to show that $b\leq a_i$.  Note that
\begin{equation}\label{lp eq}
b\leq a_i+c \text{ \,and\, } b\leq a_i+c^*,
\end{equation}
as $a\leq a_i$.  Now $[a_i,1]=\{a_i,1\}$, by (\ref{(iii')}), so $a_i+c$ is $a_i$ or $1$.  If $a_i+c=a_i$, then $b\leq a_i$, by (\ref{lp eq}).  If $a_i+c=1$, then $c^*\leq a_i$, by Lemma~\ref{1 lemma}(\ref{*}), so (\ref{lp eq}) gives $b\leq a_i+c^*=a_i$.
\end{proof}

Just as in Theorem~\ref{main theorem}, the first condition in Theorem~\ref{first cor} could be replaced by `$a = (a + c) \bcdw (a + c^{*})$ for all $a,c\in S$', without loss of strength.

Let us say that
\begin{enumerate}
\item
a reduced matrix model $\langle\sbA,F\rangle$ of\/ $\,\vdash$ is $\,\vdash$\,--\emph{simple} if the interval $[F,A]$ in the $\,\vdash$\,--filter lattice of $\sbA$ has just two elements; and that

\item
$\vdash$ is \emph{semisimple} if every $\,\vdash$\,--subdirectly irreducible reduced matrix model of\/ $\,\vdash$ is $\,\vdash$\,--simple.
\end{enumerate}
Specializing Theorem~\ref{first cor}, we get an order-theoretic explanation of the following result of L\'{a}vi\v{c}ka and P\v{r}enosil.
\begin{cor}\label{2nd cor}
\textup{(\cite{LP})}
\,Let\/ $\,\vdash$ be a
protoalgebraic deductive system.
Then\/ $\,\vdash$ has an EML
iff it has an IL and
is semisimple.
\end{cor}

Somewhat more than this can be said.
Let $\boldsymbol{\Psi}=\{\Psi_n:n\in\mathbb{N}^+\}$ be as in Definition~\ref{eml def}.  We have noted that
$\boldsymbol{\Psi}$ establishes an IL for $\,\vdash$, but it is pointed out in \cite{LP} that $\boldsymbol{\Psi}$
is actually a \emph{classical} IL-sequence for $\,\vdash$ (in the sense of \cite{Raf13}), i.e., an IL-sequence with the additional property that
\[
\Gamma\cup\Psi_n(\al_1,\dots,\al_n) \textup{ \,is inconsistent in $\,\vdash$ iff\, } \Gamma\vdash\{\al_1,\dots,\al_n\}.
\]
(Conversely, the resulting notion of a \emph{classical IL} induces an EML \cite{LP}.)
A deductive system with a classical IL has a deduction-detachment theorem
\cite[p.\,401]{Raf13}, and is therefore protoalgebraic and filter-distributive.  These facts (and Corollary~\ref{weak cor}) yield a more informative variant of Corollary~\ref{2nd cor}:

\begin{thm}\label{eml char}
The following conditions on\/ $\,\vdash$ are equivalent.
\begin{enumerate}
\item\label{eml char1}
$\,\vdash$ has an EML.

\item\label{eml char4}
$\,\vdash$ is protoalgebraic, semisimple and has an IL.

\item\label{eml char2}
$\,\vdash$ is protoalgebraic and filter-distributive and, for every algebra $\sbA$\textup{,} the join semilattice of\/ compact\/ $\,\vdash$\,--filters of\/ $\sbA$ is dually pseudo-complemented and satisfies\/
$y\bcdw y^{*}=0$.
\item\label{eml char3}
$\,\vdash$ is protoalgebraic and filter-distributive and the join semilattice of\/ compact\/ $\,\vdash$\,--theories is dually pseudo-comple\-mented and satisfies\/
$y\bcdw y^*=0$.
\end{enumerate}
\end{thm}
\noindent
The equivalence of (\ref{eml char1}) and (\ref{eml char4}) was already obtained in \cite{LP}.

\section{Weak Excluded Middle Laws in Algebraizable Logics}

\begin{defn}\label{compatibility def}
{\em
A congruence relation $\theta$ on an algebra $\sbA$ is said to be {\em compatible with\/} a subset $F$
of $A$ provided that $F$ is a union of
$\theta$--classes, i.e., whenever $a\equiv_{\,\theta}b$ and $a\in F$, then $b\in F$.
}
\end{defn}

Given a quasivariety $\mathsf{K}$ (of algebras) and an algebra $\sbA$ of the same type,
the $\mathsf{K}$--{\em congruences\/} of
$\sbA$ are the congruences $\theta$ such that $\sbA/\theta\in\mathsf{K}$.
They form an algebraic closure system over $A\times A$, and hence an algebraic lattice, ordered by inclusion, in which
the compact elements are just the finitely generated $\mathsf{K}$--congruences.
At the same time, for any subset $F$ of $\sbA$, there is always a largest congruence of $\sbA$ that is compatible with $F$ \cite[Thm.\,1.5]{BP89}; it is denoted by $\Omega^\sbA F$.  \,Thus, $\langle\sbA,F\rangle$ is reduced iff $\Omega^\sbA F=\{\langle a,a\rangle:a\in A\}$.

A deductive system $\,\vdash$ is said to be {\em algebraized\/} by
$\mathsf{K}$
if, for every algebra $\sbA$, the rule $F\mapsto\Omega^\sbA F$ defines an isomorphism from the $\,\vdash$\,--filter lattice of
$\sbA$
onto the lattice of $\mathsf{K}$--congruences
of $\sbA$.
We say that $\,\vdash$ is (elementarily) {\em algebraizable\/}
if it is algebraized by {\em some\/} quasivariety $\mathsf{K}$.  In this case $\mathsf{K}$ is unique and is called the \emph{equivalent quasivariety of} $\,\vdash$.  (It comprises the algebra reducts $\sbA$ of the reduced matrix models $\langle\sbA,F\rangle$ of $\,\vdash$.)
These definitions are equivalent to the original syntactic ones; see
Blok and Pigozzi \cite{BP89}.

Every algebraizable deductive system is protoalgebraic.  In fact, a deductive system $\,\vdash$ is protoalgebraic iff, for every algebra $\sbA$, the function $F\mapsto\Omega^\sbA F$ is inclusion-preserving on the $\,\vdash$\,--filters of $\sbA$ (see \cite{Cze01,Fon16}).

When $\mathsf{K}$ is a variety and $\sbA\in\mathsf{K}$, the congruences and $\mathsf{K}$--congruences of $\sbA$ coincide, so the prefix $\mathsf{K}$-- can be dropped (and likewise the qualifier `relatively' and its signifier `R' in the definitions to follow).

An algebra $\sbA$ is said to be \emph{trivial}
if $\left|A\right|=1$.
The following result is due to Gorbunov.  It was proved first for varieties, by Koll\'{a}r \cite{Kol79}.
\begin{thm}\label{kollar}
\textup{(\cite{Gor86},\,\cite[Thm.\,2.3.16]{Gor98})}\,
A quasivariety\/ $\mathsf{K}$ has the property that
$A\times A$ is compact in the\/ $\mathsf{K}$--congruence lattice of $\sbA,$ for all $\sbA\in\mathsf{K},$ iff the nontrivial members of\/ $\mathsf{K}$ lack trivial subalgebras.
\end{thm}
\begin{defn}\label{kollar defn}
\textup{A quasivariety satisfying the conditions of Theorem~\ref{kollar} will be called a \emph{Koll\'{a}r quasivariety}.}
\end{defn}
\noindent
Further characterizations can be found in \cite{CV12}.  Thus, a quasivariety $\mathsf{K}$ that algebraizes a deductive system $\,\vdash$ is a Koll\'{a}r quasivariety iff $\mathit{Fm}$ is compact in the lattice of $\,\vdash$\,--theories
(a feature also forced by
the existence of an IL).  Except for its use of this fact, the next result is taken from \cite{Raf13}.

\begin{thm}\label{alg case}
\textup{(\cite[Thm.~3.10]{Raf13})}\,
Let\/ $\mathsf{K}$ be a quasivariety that algebraizes a deductive system\/ $\,\vdash$\textup{.}
Then the following conditions are equivalent.
\begin{enumerate}
\item
$\,\vdash$ has an inconsistency lemma.

\item
For every algebra $\sbA$\textup{,} the join semilattice of compact\/
$\mathsf{K}$--congruences of $\sbA$ is dually pseudo-complemented.

\item
For every\/ $\sbA\in\mathsf{K}$\textup{,} the join semilattice of compact\/ $\mathsf{K}$--congruences of\/ $\sbA$
is dually pseudo-complemented.
\setcounter{newexmp}{\value{enumi}}
\end{enumerate}
In this case, the nontrivial members of\/ $\mathsf{K}$ lack trivial subalgebras.
If\/ $\mathsf{K}$ is a variety, then the numbered conditions are equivalent to
\begin{enumerate}
\setcounter{enumi}{\value{newexmp}}
\item
For every\/ $\sbA\in\mathsf{K}$\textup{,} the join semilattice of compact congruences of\/ $\sbA$
is dually pseudo-complemented.
\end{enumerate}
\end{thm}

Given a quasivariety $\mathsf{K}$, we say that an algebra $\sbA\in\mathsf{K}$ is
\emph{relatively subdirectly irreducible} (RSI),
or {\em relatively finitely subdirectly irreducible\/} (RFSI),
or \emph{relatively simple} (RS)
if, in the lattice of $\mathsf{K}$--congruences of $\sbA$,
the
relation $
\{\langle a,a\rangle :
a\in A\}$
is completely meet-irreducible, or meet-irreducible,
or a co-atom,
respectively.
The class of all RSI [RFSI; RS]
algebras in $\mathsf{K}$
is denoted by
$\mathsf{K}_\textup{RSI\/}$ [$\mathsf{K}_\textup{RFSI\/}$; $\mathsf{K}_\textup{RS\/}$].
Thus, $
\mathsf{K}_\textup{RS\/}\subseteq
\mathsf{K}_\textup{RSI\/}\subseteq \mathsf{K}_\textup{RFSI\/}$, and $\mathsf{K}_\textup{RSI\/}$ consists of nontrivial algebras.
Every member of $\mathsf{K}$ is isomorphic to a subdirect product of members of
$\mathsf{K}_\textup{RSI\/}$ \cite[Thm.~1.1]{Pig88}.
If every RSI member of $\mathsf{K}$ is relatively simple, then $\mathsf{K}$ is said to be \emph{relatively semisimple}.

When $\mathsf{K}$ is the equivalent quasivariety of $\,\vdash$, then $\mathsf{K}_\textup{RSI\/}$ [$\mathsf{K}_\textup{RS\/}$; $\mathsf{K}_\textup{RFSI\/}$]
comprises the algebra reducts of the reduced matrix models $\langle\sbA,F\rangle$ of $\,\vdash$ that are $\,\vdash$\,--subdirectly irreducible [that are $\,\vdash$\,--simple; for which $F$ is meet-irreducible in the $\,\vdash$\,--filter lattice of $\sbA$].
In this case, $\,\vdash$ is semisimple in the sense of Section~\ref{excluded middle laws} iff $\mathsf{K}$ is relatively semisimple.

For algebraizable logics, our main result about the WEML (Theorem~\ref{prsi greatest}) therefore takes the following form.

\begin{thm}\label{rsi greatest}
Let\/ $\mathsf{K}$ be a quasivariety algebraizing a deductive system\/ $\,\vdash$ with an inconsistency lemma.  Then the following conditions are equivalent.
\begin{enumerate}
\item\label{rsi greatest1}
$\,\vdash$ has a WEML.

\item\label{rsi greatest2}
Every nontrivial algebra in\/ $\mathsf{K}_\textup{RFSI}$ has a greatest proper\/ $\mathsf{K}$--congruence (i.e., it has a greatest proper\/ $\,\vdash$\,--filter).

\item\label{rsi greatest3}
Every algebra in\/ $\mathsf{K}_\textup{RSI}$ has a greatest proper\/ $\mathsf{K}$--congruence.
\end{enumerate}
\end{thm}

\begin{cor}
A super-intuitionistic logic has a WEML iff it extends\/ $\mathbf{KC}$ \textup{(}i.e., its theorems include $\neg v\vee\neg\neg v).$
\end{cor}

\begin{cor}\label{cat eq}
Let\/ $\mathsf{K}_1$ and\/ $\mathsf{K}_2$ be categorically equivalent quasivarieties that algebraize deductive systems\/ $\,\vdash_1$ and\/ $\,\vdash_2,$ respectively.  If\/ $\,\vdash_1$ has a WEML, then so does\/ $\,\vdash_2.$
\end{cor}
\begin{proof}
A category equivalence functor $F$ from $\mathsf{K}_1$ to $\mathsf{K}_2$ induces an isomorphism from the $\mathsf{K}_1$--congruence lattice of each $\sbA\in\mathsf{K}_1$ onto the $\mathsf{K}_2$--congruence lattice of its image $F(\sbA)\in\mathsf{K}_2$.  Also, a lattice isomorphism between complete lattices restricts to an isomorphism between their join semilattices of compact elements.  The result therefore follows from Theorems~\ref{alg case} and \ref{rsi greatest}.
\end{proof}

Corollary~\ref{cat eq} applies equally to EMLs (cf.\ \cite[Cor.\,5.7]{Raf13}), because they amount to classical ILs.
When $\,\vdash$ is algebraized by a quasivariety $\mathsf{K}$, then
it has a classical IL iff $\mathsf{K}$ is a relatively filtral Koll\'{a}r quasivariety \cite{CR17}.
For the definition of relative filtrality (and some alternative characterizations),
see \cite{CR17} and its references.

\begin{exmp}\label{product logic}
\textup{The {\em product logic\/} $\boldsymbol{\Pi}$ of \cite{Haj98} is algebraizable and
has an inconsistency lemma, with
$\{\{\neg(v_1\odot\,\dots\,\odot v_n)\}:n\in\mathbb{N}^+\}$
as IL-sequence.  It therefore has a WEML, by Theorem~\ref{rsi greatest}, because the subdirectly irreducible members of its equivalent variety are totally ordered and have least elements, and their deductive filters are upward-closed.  (In such an algebra, the union of the proper deductive filters is the largest proper deductive filter.
Alternatively, one can argue syntactically from the theoremhood of $\neg v\vee\neg\neg v$ in $\boldsymbol{\Pi}$.)
On the other hand, $\boldsymbol{\Pi}$ does not have an EML, as it lacks a classical IL \cite[p.\,404]{Raf13} (equivalently, it is not semisimple).}
\end{exmp}

\section{Normal Modal Logics}

Recall that a \emph{modal formula} is a formula in the signature of classical propositional logic, expanded by a unary connective $\Box$, where $\neg\Box\neg\al$ is abbreviated as $\Diamond\al$.  We define $\Box^0\al=\al$ and $\Box^{n+1}\al=\Box\Box^n\al$ for $n\in\omega=\{0,1,2,\dots\}$, and similarly for $\Diamond$.  Moreover, for $n\in\omega$, we define
\begin{align*}
\boxplus^n\varphi \,=\, \varphi \,\wedge\, \Box\varphi
\,\wedge\, {\dots} \,\wedge\,
\Box^{n}\varphi;\\
\diamondplus^n\varphi \,=\, \varphi \,\vee\, \Diamond\varphi \,\vee\, {\dots} \,\vee\, \Diamond^n\varphi.
\end{align*}
A \emph{normal modal logic} $\mathbf{L}$ is traditionally identified with a special set of modal formulas, rather than a consequence relation.  More precisely, $\mathbf{L}$ must include all classical tautologies and Kripke's distribution axiom
\[
\Box(v_1\to v_2)\to(\Box v_1\to \Box v_2),
\]
and $\mathbf{L}$ must be closed under substitution, under modus ponens, and under the connective $\Box$ (i.e., under the rule of necessitation).  We denote by $\,\vdash_\mathbf{L}$ the \emph{global consequence relation of} $\mathbf{L}$ (see \cite{Kra07}).  Thus, the elements of $\mathbf{L}$ are exactly the theorems of $\,\vdash_\mathbf{L}$, provided that modus ponens and necessitation serve as the inference rules of $\,\vdash_\mathbf{L}$.

For $m,n\in\omega$, we then have ${\Diamond^n}v\leftrightarrow\neg{\Box^n}\neg v,\, {\diamondplus^n}v\leftrightarrow\neg{\boxplus^n}\neg v\in\mathbf{L}$, and if $m\leq n$, then ${\boxplus^n}v\to{\boxplus^m}v,\,{\diamondplus^m}v\to{\diamondplus^n}v\in\mathbf{L}$.
For each normal modal logic $\mathbf{L}$, the system
$\,\vdash_\mathbf{L}$ has the following \emph{local deduction-detachment theorem} (LDDT):
\[
\textup{$\Gamma\cup\{\al\}\vdash_\mathbf{L}\be$ iff there exists $n\in\omega$ such that
$\Gamma\vdash_\mathbf{L}\boxplus^n\al\to\be$.}
\]

We characterize below the normal modal logics $\mathbf{L}$ for which $\,\vdash_\mathbf{L}$ has an IL or a WEML.  (The semisimple systems of this kind are already understood \cite{KK06,LP}, so the case of an EML requires no further attention.)

\begin{thm}\label{modal il}
Let\/ $\mathbf{L}$ be a normal modal logic.  Then
$\,\vdash_\mathbf{L}$ has an inconsistency lemma iff
there exists $n\in\omega$ such that\/ $\,\vdash_\mathbf{L} {\boxplus^n v} \to {\diamondplus^n}{\boxplus^{n+1}v}$
(i.e., ${\boxplus^n \al} \to {\diamondplus^n}{\boxplus^{n+1}\al}\in\mathbf{L}$ for all modal formulas $\al)$.
\end{thm}
\begin{proof}
Observe first that, for any $\Gamma\cup\{\al_1,\dots,\al_k\}\subseteq\mathit{Fm}$,
\begin{align}\label{modal 0}
\begin{array}{l}
\textup{$\Gamma\cup\{\al_1,\dots,\al_k\}$
is inconsistent
 in $\,\vdash_\mathbf{L}$}\\[0.2pc] \textup{iff
\,$\Gamma\cup\{\al_1\wedge{\dots}\wedge\al_k\}\vdash_\mathbf{L}\bot$,}\\[0.2pc]
\textup{iff
\,${\Gamma\vdash_\mathbf{L}\neg{\boxplus}^m(\al_1\wedge{\dots}\wedge\al_k)}$ for some $m\in\omega$}
\end{array}
\end{align}
(in view of the LDDT).

($\Rig$)\, Let $
\{\Psi_n:n\in\mathbb{N}^+\}$ be an elementary IL-sequence for $\,\vdash_\mathbf{L}$.
Because ${\Psi_1(v)\cup\{v\}}$
is inconsistent in $\,\vdash_\mathbf{L}$, (\ref{modal 0}) shows that we can choose $n\in\omega$ with
${\Psi_1(v)\vdash_\mathbf{L}\neg{\boxplus^n}v}$.
Also, $\{\neg{\boxplus^{n+1}}v,v\}$ is inconsistent in $\,\vdash_\mathbf{L}$ (owing to necessitation), so $\neg{\boxplus^{n+1}}v\vdash_\mathbf{L}\Psi_1(v)$, by the IL, whence
$\neg{\boxplus^{n+1}}v\vdash_\mathbf{L}\neg{\boxplus^n}v$.  Thus, by the LDDT, there exists $m\in\omega$ such that
$\,\vdash_\mathbf{L}{\boxplus^m}\neg{\boxplus^{n+1}}v \to\neg{\boxplus^n} v$, i.e.,
$\,\vdash_\mathbf{L}{\boxplus^n} v \to \neg {\boxplus^m}\neg{\boxplus^{n+1}}v$, i.e.,
\begin{equation}\label{modal 1}
\,\vdash_\mathbf{L}{\boxplus^n} v \to {\diamondplus^m}{\boxplus^{n+1}}v.
\end{equation}
If $m\leq n$ then
$\,\vdash_\mathbf{L}{\diamondplus^m}{\boxplus^{n+1}}v\to {\diamondplus^n}{\boxplus^{n+1}}v$, whence
$\,\vdash_\mathbf{L} {\boxplus^n} v \to {\diamondplus^n}{\boxplus^{n+1}}v$, in view of (\ref{modal 1}).  And if $n<m$ then, substituting ${\boxplus^{m-n}}v$ for $v$ in (\ref{modal 1}),
we obtain
$\,\vdash_\mathbf{L} {\boxplus^m} v \to {\diamondplus^m}{\boxplus^{m+1}}v$.

($\Leftarrow$)\,
Let $n\in\omega$ be such that
\begin{equation}\label{modal 4}
\,\vdash_\mathbf{L} {\boxplus^n} v \to {\diamondplus^n}{\boxplus^{n+1}}v.
\end{equation}
For each positive integer $k$, let $\Psi_k=\{\neg{\boxplus^n}(v_1\wedge{\dots}\wedge v_k)\}$.  We shall show that $\{\Psi_k:k\in\mathbb{N}^+\}$ is an IL-sequence for $\,\vdash_\mathbf{L}$.  Let $k\in\mathbb{N}^+$ and suppose $\Gamma\cup\{\al_1,\dots,\al_k\}\subseteq\mathit{Fm}$.

If $\Gamma\vdash_\mathbf{L}\Psi_k(\al_1,\dots,\al_k)$, then $\Gamma\cup\{\al_1,\dots,\al_k\}$ is inconsistent in $\vdash_\mathbf{L}$, by (\ref{modal 0}).  We need to prove the converse, so we may assume (again by (\ref{modal 0})) that there exists $m\in\omega$ with
\begin{equation}\label{modal 2}
\Gamma\vdash_\mathbf{L}\neg{\boxplus^m}(\al_1\wedge{\dots}\wedge\al_k).
\end{equation}
If $m\leq n$ then $\neg{\boxplus^m}(\al_1\wedge{\dots}\wedge\al_k)\vdash_\mathbf{L}\Psi_k(\al_1,\dots,\al_k)$, in which case, by (\ref{modal 2}), ${\Gamma\vdash_\mathbf{L}\Psi_k(\al_1,\dots,\al_k)}$, as required.  We may therefore assume that ${n<m}$.
By necessitation, $\neg{\boxplus^{n+1}}v \vdash_\mathbf{L} {\boxplus^n} \neg{\boxplus^{n+1}}v$, but
${\boxplus^n} \neg{\boxplus^{n+1}}v  \vdash_\mathbf{L} \neg{\diamondplus^n}{\boxplus^{n+1}}v$, so
\begin{equation}\label{modal 3}
\neg{\boxplus^{n+1}}v \vdash_\mathbf{L} \neg{\diamondplus^n}{\boxplus^{n+1}}v.
\end{equation}
By (\ref{modal 4}) and contraposition,
$\,\vdash_\mathbf{L} \neg{\diamondplus^n}{\boxplus^{n+1}}v \to \neg {\boxplus^n}v$, so by (\ref{modal 3}),
\[
\neg{\boxplus^{n+1}}v \vdash_\mathbf{L} \neg {\boxplus^n}v.
\]
Therefore, because $n<m$, the substitution-invariance and transitivity of $\,\vdash_\mathbf{L}$ yield $\neg{\boxplus^m}v \vdash_\mathbf{L} \neg {\boxplus^n}v$.
This, with (\ref{modal 2}), gives $\Gamma\vdash_\mathbf{L}\Psi_k(\al_1,\dots,\al_n)$, as required.
\end{proof}

\begin{thm}\label{modal weml}
Let\/ $\mathbf{L}$ be a normal modal logic.  Then\/ $\,\vdash_\mathbf{L}$ has a WEML iff there exists $n\in\omega$ such that\/
\[
\textup{$\vdash_\mathbf{L} {\boxplus^n v} \to {\diamondplus^n}{\boxplus^{n+1}v}$\, and\/
$\,\vdash_\mathbf{L}{\boxplus^m}\neg{\boxplus^n}v\,\vee\,{\boxplus^m}\neg{\boxplus^n}\neg{\boxplus^n}v,$ for all\/ $m\in\omega.$}
\]
\end{thm}
\begin{proof}
By Theorem~\ref{modal il} and its proof, we may assume that
\[
\vdash_\mathbf{L} {\boxplus^n v} \to {\diamondplus^n}{\boxplus^{n+1}v}
\]
for some (fixed) $n\in\omega$, whence
$\boldsymbol{\Psi}=\{\Psi_k:k\in\mathbb{N}^+\}$ is an IL-sequence for $\,\vdash_\mathbf{L}$,
where $\Psi_k=\{\neg{\boxplus^n}(v_1\wedge{\dots}\wedge v_k)\}$ for each $k$.

($\Leftarrow$)\,
To prove that $\boldsymbol{\Psi}$ establishes a WEML for $\,\vdash_\mathbf{L}$, suppose
\[
\textup{$\Gamma\cup\{\neg{\boxplus^n}\al\}\vdash_\mathbf{L}\be$ \ and \ $\Gamma\cup\{\neg{\boxplus^n}\neg{\boxplus^n}\al\}\vdash_\mathbf{L}\be$,}
\]
where $\al$ is $\al_1\wedge{\dots}\wedge\al_k$ ($k\in\mathbb{N}^+$).
By the LDDT, there exists $m\in\omega$ with
\[
\textup{$\Gamma\vdash_\mathbf{L}{\boxplus^m}\neg{\boxplus^n}\al\to\be$ \ and \ $\Gamma\vdash_\mathbf{L}{\boxplus^m}\neg{\boxplus^n}\neg{\boxplus^n}\al\to\be$.}
\]
With the help of a classical tautology, we obtain
\[
\Gamma\vdash_\mathbf{L}({\boxplus^m}\neg{\boxplus^n}\al\,\vee\,{\boxplus^m}\neg{\boxplus^n}\neg{\boxplus^n}\al)\to\be.
\]
Then, by the given assumption and modus ponens, $\Gamma\vdash\be$, as required.

($\Rig$)\, Let $m\in\omega$.  By necessitation, $\neg{\boxplus^n}v\vdash_\mathbf{L} {\boxplus^m}\neg{\boxplus^n}v$, and so
\[
\neg{\boxplus^n}v \vdash_\mathbf{L} {\boxplus^m}\neg{\boxplus^n}v\,\vee\, {\boxplus^m}\neg{\boxplus^n}\neg{\boxplus^n}v,
\]
i.e., $\Psi_1(v) \vdash_\mathbf{L} {\boxplus^m}\neg{\boxplus^n}v\,\vee\, {\boxplus^m}\neg{\boxplus^n}\neg{\boxplus^n}v$.
Similarly,
\[
\Psi_1\Psi_1(v) \vdash_\mathbf{L} {\boxplus^m}\neg{\boxplus^n}v\,\vee\, {\boxplus^m}\neg{\boxplus^n}\neg{\boxplus^n}v,
\]
so by the WEML, $\,\vdash_\mathbf{L} {\boxplus^m}\neg{\boxplus^n}v\,\vee\, {\boxplus^m}\neg{\boxplus^n}\neg{\boxplus^n}v$.
\end{proof}

\section{Extensions of $\mathbf{S4}$}

Recall that $\mathbf{S4}$ is the smallest normal modal logic $\mathbf{L}$ such that
\[
{\Box v\to v},\,{\Box v \to \Box\Box v}\in\mathbf{L}
\]
(equivalently, $v\to\Diamond v,\,\Diamond\Diamond v\to \Diamond v\in\mathbf{L}$),
and that $\textup{Next}(\mathbf{S4})$ is the lattice of normal modal logics containing $\mathbf{S4}$.  These logics prove ${\boxplus^m}v\leftrightarrow\Box v$ and ${\diamondplus^m}v\leftrightarrow\Diamond v$ for all $m\in\mathbb{N}^+$.

If $\mathbf{L}\in\textup{Next}(\mathbf{S4})$, then since $\Box v\to \Diamond\Box v\in\mathbf{L}$, we have
${\boxplus^1}v\to{\diamondplus^1}{\boxplus^2}v\in\mathbf{L}$.  Therefore,
Theorem~\ref{modal il} and its proof yield the following.
\begin{exmp}\label{s4 il}
\textup{For each\/ $\mathbf{L}\in\textup{Next}(\mathbf{S4}),$ the global consequence relation\/ $\,\vdash_\mathbf{L}$ has an inconsistency lemma,
with\/ $n=1$ in Theorem~\textup{\ref{modal il},}
and an IL-sequence\/ $\{\Psi_k:k\in\mathbb{N}^+\}$ for $\,\vdash_\mathbf{L}$ is given by\/
$\Psi_k=\{\neg\Box(v_1\wedge{\dots}\wedge v_k)\}.$}
\end{exmp}

We shall show that, for each $\mathbf{L}\in\textup{Next}(\mathbf{S4})$, the system $\,\vdash_\mathbf{L}$ has a WEML iff its theorems include the so-called \emph{convergence axiom} $\Diamond \Box v \to \Box \Diamond v$.  This formula is validated by a Kripke frame $\sbX=\langle X,R\rangle$ iff $\sbX$ is \emph{principally up-directed} in the following sense:
\begin{equation*}
\begin{array}{l}
\textup{for any $x,y,z \in X$ such that $xRy$ and $xRz$,}\\[0.1pc]
\textup{there exists $w \in X$ such that $yRw$ and $zRw$.}
\end{array}
\end{equation*}
The extension of $\mathbf{S4}$ by the convergence axiom is known as $\mathbf{S4.2}$.  Thus, $\mathbf{S4.2}$
is the normal modal logic induced by the class of Kripke frames that are reflexive, transitive and principally up-directed.  On the other hand, $\mathbf{KC}$ is the super-intuitionistic logic induced by the class of
principally
up-directed posets.
It was proved in \cite{DL59} that $\mathbf{S4.2}$ is the least modal companion of $\mathbf{KC}$ (see \cite{CZ92,MR74,Ryb97} for the general notion of a modal companion and further examples).

\begin{thm}\label{s4 weml}
Let\/ $\mathbf{L}\in\textup{Next}(\mathbf{S4}).$  Then the global consequence relation\/ $\,\vdash_\mathbf{L}$ has a WEML iff\/ $\,\vdash_\mathbf{L}\Diamond \Box v \to \Box \Diamond v$ \textup{(}i.e., $\mathbf{L}$ extends\/ $\mathbf{S4.2}).$
\end{thm}
\begin{proof}
($\Rig$)\,
Suppose $\,\vdash_\mathbf{L}$ has a WEML.  By Theorem~\ref{modal weml} and Example~\ref{s4 il},
\[
\textup{$\,\vdash_\mathbf{L} \,{\boxplus^m} \lnot \Box v \,\,\lor\,\, {\boxplus^m} \lnot \Box \lnot \Box v$, \,\,for all $m \in \omega$.}
\]
For $m=1$, this gives
$\,\vdash_\mathbf{L}\Box \lnot \Box v \,\lor\, \Box \lnot \Box \lnot \Box v$,
i.e., $\,\vdash_\mathbf{L}\lnot \Diamond \Box v \,\lor\, \Box \Diamond \Box v$, i.e.,
\begin{equation}\label{s4 1}
\,\vdash_\mathbf{L}\Diamond \Box v \to \Box \Diamond \Box v.
\end{equation}
We also have $v_2\to v_3\vdash_\mathbf{L}\Box\Diamond v_2\to\Box\Diamond v_3$ (in any normal modal logic), so from
$\,\vdash_\mathbf{L}\Box v \to v$,
we may infer
$\,\vdash_\mathbf{L}\Box \Diamond \Box v \to \Box \Diamond v$.  This, with (\ref{s4 1}), shows that $\mathbf{L}$ includes
the convergence axiom.

($\Leftarrow$)\,
Suppose $\,\vdash_\mathbf{L}\Diamond \Box v \to \Box \Diamond v$.
Substituting $\Box v$ for $v$, we obtain
\[
\,\vdash_\mathbf{L}\Diamond \Box \Box v \to \Box \Diamond \Box v.
\]
We now use repeatedly, without comment, the fact that $\mathbf{L}\in \textup{Next}(\mathbf{S4})$.  The formulas $\Diamond \Box \Box v$ and $\Diamond \Box v$ are logically equivalent over $\mathbf{L}$, whence
\[
\textup{$\,\vdash_\mathbf{L}\Diamond \Box v \to \Box \Diamond \Box v$, \,\,i.e.,\,\,
$\,\vdash_\mathbf{L}\Box \lnot \Box v \,\lor\, \Box \lnot \Box \lnot \Box v$.}
\]
This implies that
$\,\vdash_\mathbf{L}{\boxplus^m}(\lnot \Box v) \,\lor\, {\boxplus^m}(\lnot \Box \lnot \Box v)$
for every $m \geq 1$.  Furthermore,
${\boxplus}^0(\lnot \Box v) \,\lor\, {\boxplus^0}(\lnot \Box \lnot \Box v)$ is logically equivalent, over $\mathbf{L}$, to $\Box v \to \Diamond \Box v$, which belongs to $\mathbf{L}$.
Thus, for all $m \in \omega$, we have
\[
\,\vdash_\mathbf{L}{\boxplus^m}(\lnot \Box v) \,\lor\, {\boxplus^m}(\lnot \Box \lnot \Box v).
\]
The formula $\Box v$ is logically equivalent, over $\mathbf{L}$, to ${\boxplus^1}v$.  Therefore,
\[
\textup{$\,\vdash_\mathbf{L}{\boxplus^m}(\lnot {\boxplus^1}v) \,\lor\, {\boxplus^m}(\lnot {\boxplus^1} \lnot {\boxplus^1}v)$ \,for all $m\in\omega$.} \]
This, with Theorem~\ref{modal weml} and Example~\ref{s4 il}, shows that
$\,\vdash_\mathbf{L}$ has a WEML.
\end{proof}

It is now easy to construct logics $\mathbf{L}\in\textup{Next}(\mathbf{S4})$ for which $\,\vdash_\mathbf{L}$ has a WEML but lacks an EML. Indeed, let $\sbX$ be any Kripke frame that is reflexive, transitive and principally up-directed.
Let
$\mathbf{L}$ be the normal modal logic
induced by $\sbX$.
Then $\mathbf{L}$
extends $\mathbf{S4.2}$ and
$\,\vdash_\mathbf{L}$
has a WEML, by Theorem~\ref{s4 weml}.  If we assume, moreover, that $\sbX$ is rooted and contains points $x,y$ such that $x{\not\mathrel{R}} y$ or $y{\not\mathrel{R}}x$, then the complex algebra of $\sbX$ is subdirectly irreducible but not simple. Consequently, $\,\vdash_\mathbf{L}$
lacks a EML, by Theorem~\ref{eml char}.  A concrete example is the case where $\sbX$ is a two-element chain, viewed as a poset.

\section{Relevance Logics}\label{dmm section}

Dunn \cite{Dun66,MDL74} showed in 1966 that
the variety $\mathsf{DMM}$ of De Morgan monoids algebraizes
the principal relevance logic $\mathbf{R^t}$.
(More exactly, it algebraizes the deducibility relation $\,\vdash_{\mathbf{R^t}}$ of the formal system $\mathbf{R^t}$ from \cite{AB75}, but we shall often abbreviate that relation as $\mathbf{R^t}$.)
Consequently, the subvarieties of $\mathsf{DMM}$ and the axiomatic extensions of $\mathbf{R^t}$ form anti-isomorphic lattices.

We shall show that an axiomatic extension of
$\mathbf{R^t}$ has an IL iff it is algebraized by a Koll\'{a}r variety of De Morgan monoids, and that in this case it also has a WEML.  Some characterizations of the Koll\'{a}r subvarieties of $\mathsf{DMM}$ will be provided.  Where known structural features of De Morgan monoids are mentioned below without citation, their sources are given in the recent papers \cite{MRW19,MRW20b,MRW20}.

\begin{defn}\label{dmm def}
\textup{A
\emph{De Morgan monoid} is an algebra
$\sbA=\langle A;\bcdw,\wedge,\vee,\neg,e\rangle$
comprising a distributive lattice $\langle A;\wedge,\vee\rangle$,
a commutative monoid $\langle A;\bcdw,e\rangle$ that is \emph{square-increasing}
(i.e., $\sbA$
satisfies $x\leqslant x^2\seteq x\bcdw x$),
and a function $\neg\colon A\mrig A$,
called an {\em involution},
such that $\sbA$ satisfies $\neg\neg x=
x$ and
\begin{equation*}
x\bcdw y\leqslant z\;\Longleftrightarrow\;x\bcdw\neg z\leqslant\neg y.
\end{equation*}
Here, $\al\leqslant\be$ abbreviates $\al=
\al\wedge\be$.
We refer to $\bcdw$ as {\em fusion}, and we define
\[
\textup{$f=\neg e$ \,\,and\,\, $x\to y=\neg(x\bcdw\neg y)$ \,\,and\,\, $x\leftrightarrow y=(x\to y)\wedge (y\to x)$.}
\]}
\end{defn}
It follows that $\neg$ is an anti-automorphism of $\langle A;\wedge,\vee\rangle$ (so De Morgan's laws hold), and that
$\sbA$ satisfies
the \emph{law of residuation\/}:
\[
x\bcdw y\leqslant z\;\Longleftrightarrow\;y\leqslant x\rig z.
\]
In particular, $\sbA$ satisfies
\begin{equation}\label{order}
x\leqslant z\;\Longleftrightarrow\;e\leqslant x\rig z,
\end{equation}
as well as $e\to x=x$.  And (\ref{order}) shows that $e$ is not the least element of $\sbA$, unless $\sbA$ is trivial.

It turns out that the
$\mathbf{R^t}$--filters of a De Morgan monoid
$\sbA$ are just the lattice-filters $F$ of
$\langle A;\wedge,\vee\rangle$ such that $e\in F$.  The smallest of these is therefore $[e)\seteq\{a\in A:e\leqslant a\}$.
An $\mathbf{R^t}$--filter of $\sbA$ is closed under fusion, owing to the square-increasing law.
The natural lattice isomorphism from
$\mathbf{R^t}$--filters to congruences of $\sbA$, and its inverse, are given by
\begin{eqnarray*}
& F\,\mapsto\,\Omega^\sbA F=\{\langle a,b\rangle\in A\times A: a\leftrightarrow b
\in F\}
;\\
& \theta\,\mapsto\,\{a\in A:
a\wedge e\equiv_\theta e
\}.
\end{eqnarray*}
The deductive system $\,\vdash_\mathbf{R^t}$ can in fact be characterized as the consequence relation of the class of matrices $\{\langle\sbA,[e)\rangle:\sbA=\langle A;\bcdw,\wedge,\vee,\neg,e\rangle\in\mathsf{DMM}\}$.  Partly for this reason, we shall not discuss syntactic postulates for $\mathbf{R^t}$ here.

The reader should recall Definition~\ref{kollar defn} at this point.

\begin{lem}\label{dmm lem 1}
Let\/ $\mathsf{K}$ be a Koll\'{a}r variety of De Morgan monoids, with ${\sbA\in\mathsf{K}.}$  Then $\sbA$ is bounded, i.e., the lattice\/ $\langle A;\wedge,\vee\rangle$ has a least and a greatest element.
\end{lem}
\begin{proof}
As $\mathsf{K}$ is a Koll\'{a}r variety, and since there is a lattice isomorphism between the $\mathbf{R^t}$--filters and the congruences of $\sbA$, the total $\mathbf{R^t}$--filter $A$ of $\sbA$ is compact in the $\mathbf{R^t}$--filter lattice of $\sbA$.
Whenever $e\geqslant a\in A$, then $[a)=\{b\in A:a\leqslant b\}$ is an $\mathbf{R^t}$--filter of $\sbA$, and $A$ is clearly the join of $\{[a):e\geqslant a\in A\}$.  By compactness, therefore, $A$ is already the join of $\{[a_i):i=1,\dots,n\}$ for some $a_1,\dots,a_n\leqslant e$, with $n\in\mathbb{N}^+$.  Thus, $A=[a)$, where $a\seteq a_1\wedge{\dots}\wedge a_n$, i.e., $a$ is the least element of $\sbA$, whence $\neg a$ is the greatest element.
\end{proof}

Suppose $\bot,\top$ are, respectively, the least and the greatest element of a De Morgan monoid $\sbA$.  Then $a\bcdw\bot=\bot$ for all $a\in A$, and the following conditions are equivalent:
\begin{enumerate}
\item
$a\bcdw\top=\top$ for all $a\in A\backslash\{\bot\}$;
\item
$a\to\bot=\bot$ for all $a\in A\backslash\{\bot\}$.
\end{enumerate}
When these conditions hold, we say that $\sbA$ is \emph{rigorously compact}.  In that case, no proper congruence of $\sbA$ identifies $\bot$ with an element of $A\backslash\{\bot\}$ (see \cite[Lem.\,2.1(i)]{MRW20b} or \cite[Prop.\,6.2(i)]{OR07}).  Moreover, every bounded FSI De Morgan monoid is rigorously compact (see \cite[Thm.\,5.3]{MRW19}, which has an antecedent in \cite[Thm.\,3]{Mey86}).

\begin{lem}\label{dmm lem 2}
Let\/ $\bot$ be the least element of a nontrivial FSI De Morgan monoid $\sbA.$  Then\/ $\bot$ is meet-irreducible in the sublattice $(e]=\{a\in A:a\leqslant e\}$ of\/ $\sbA.$  Consequently, $\sbA$ has a largest proper congruence.
\end{lem}
\begin{proof}
Note that $\bot<e$, as $\sbA$ is nontrivial.
By the above remarks, $\sbA$ is rigorously compact, and
$\{\bot\}$ is an equivalence class of every proper congruence of $\sbA$. Suppose, with a view to contradiction, that $a\wedge b=\bot$,
where $\bot < a, b < e$.  Let $\theta$ be the congruence $\Omega^\sbA[a)$ of $\sbA$.
As $e\leqslant a\to e$, we have $a\leftrightarrow e=(a\to e)\wedge a=a$, so $a\equiv_\theta e$.
Therefore, $\bot=a\wedge b\equiv_\theta e\wedge b=b$, so the $\theta$--class of $\bot$ is not a singleton.  Consequently, $\theta=A\times A$.  In particular, $a\equiv_\theta \bot$, i.e., $a\leqslant a\leftrightarrow\bot\leqslant a\to\bot=\bot$, i.e., $a=\bot$, a contradiction.

This confirms that $\bot$ is meet-irreducible in $(e]$, so
$(e] \backslash \{\bot\}$ is a proper filter of the lattice $(e]$. The upward closure in $\sbA$ of $(e] \backslash \{\bot\}$ is therefore the largest proper $\mathbf{R^t}$--filter of $\sbA$, and so $\sbA$ has a greatest proper congruence.
\end{proof}

The following \emph{deduction-detachment theorem} (DDT) applies to $\mathbf{R^t}$:
\begin{equation}\label{ddt for rt}
\textup{$\Gamma\cup\{\al\}\vdash_\mathbf{R^t}\be$ \,iff\, $\Gamma\vdash_\mathbf{R^t} (\al\wedge e)\to\be$.}
\end{equation}
A deductive system $\,\vdash$ with a DDT (in the general sense of \cite{BP}) has an inconsistency lemma iff $\mathit{Fm}$ is compact in the lattice of $\,\vdash$\,--theories; see \cite[Cor.\,3.9]{Raf13}.  This demand amounts, when a variety $\mathsf{K}$ algebraizes $\,\vdash$, to the requirement that $\mathsf{K}$ be a Koll\'{a}r variety.  Putting this together with Lemmas~\ref{dmm lem 1} and \ref{dmm lem 2} and Theorem~\ref{rsi greatest}, we obtain the following.
\begin{thm}\label{r weml}
An axiomatic extension of\/ $\mathbf{R^t}$ has an inconsistency lemma iff it has a WEML, iff it is algebraized by a Koll\'{a}r variety of De Morgan monoids.
\end{thm}
The following remarks illuminate the content of Theorem~\ref{r weml}.

A Boolean algebra may be regarded as a De Morgan monoid in which $\bcdw$ duplicates $\wedge$.  More generally, a De Morgan monoid $\sbA$ is \emph{idempotent} (in the sense that $a^2=a$ for all $a\in A$) iff it satisfies $f\leqslant e$\,; for a proof, see \cite[Thm.\,3.3]{MRW19}.  An \emph{odd Sugihara monoid} is a De Morgan monoid in which $f=e$.
We depict below the two-element Boolean algebra $\mathbf{2}$,
the three-element odd Sugihara monoid $\sbS_3$, and two
four-element
De Morgan monoids, $\sbC_4$ and $\sbD_4$.
In each case, the labeled Hasse diagram determines the structure.  Note that $\neg f^2$ abbreviates $\neg(f^2)$.

{\tiny

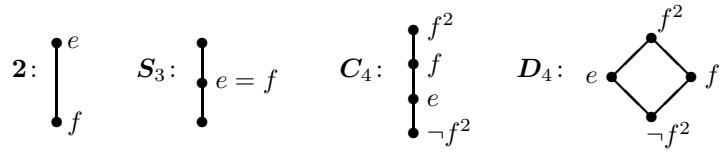
\thicklines
\begin{center}
\begin{picture}(80,60)(-28,51)

\put(-105,63){\line(0,1){30}}
\put(-105,63){\circle*{4}}
\put(-105,93){\circle*{4}}

\put(-101,91){\small $e$}
\put(-101,60){\small $f$}

\put(-122,80){\small $\mathbf{{2}\colon}$}

%
%

\put(-50,78){\circle*{4}}
\put(-50,63){\line(0,1){30}}
\put(-50,63){\circle*{4}}
\put(-50,93){\circle*{4}}

\put(-45,76){\small ${e}=f$}

\put(-75,80){\small ${\sbS_3}\colon$}

%
%

\put(30,59){\circle*{4}}
\put(30,59){\line(0,1){39}}
\put(30,72){\circle*{4}}
\put(30,85){\circle*{4}}
\put(30,98){\circle*{4}}

\put(35,96){\small ${f^2}$}
\put(35,82){\small $f$}
\put(35,69){\small ${e}$}
\put(35,56){\small $\neg f^2$}

\put(2,80){\small ${\sbC_4}\colon$}

%
%

\put(120,65){\circle*{4}}
\put(135,80){\line(-1,-1){15}}
\put(135,80){\circle*{4}}
\put(105,80){\line(1,-1){15}}
\put(105,80){\circle*{4}}
\put(105,80){\line(1,1){15}}
\put(120,95){\circle*{4}}
\put(135,80){\line(-1,1){15}}

\put(122,99){\small ${f^2}$}
\put(95,78){\small ${e}$}
\put(140,78){\small $f$}
\put(118,55){\small $\neg f^2$}

\put(69,80){\small ${\sbD_4}\colon$}


\end{picture}\nopagebreak
\end{center}

}

\noindent
As it happens, the varieties generated, respectively, by these four algebras are exactly the minimal (nontrivial) subvarieties of $\mathsf{DMM}$ \cite[Thm.\,6.1]{MRW19}.

A quasivariety $\mathsf{K}$ of De Morgan monoids is a Koll\'{a}r quasivariety iff $\sbS_3\notin\mathsf{K}$ \cite[Thm.\,8.4(iii)]{MRW20}.  Many such non-semisimple varieties are exhibited in \cite{MRW20b}.

A De Morgan monoid $\sbA$ is said to be \emph{anti-idempotent} if it satisfies $x\leqslant f^2$ (and therefore also $\neg f^2\leqslant x$).  By \cite[Cor.\,3.6]{MRW19}, this amounts to the demand that no nontrivial idempotent algebra belongs to the variety generated by $\sbA$.  In particular, $\sbC_4$ and $\sbD_4$ have this property.

If $\mathsf{K}$ is a Koll\'{a}r variety of De Morgan monoids and $\sbA\in\mathsf{K}_\textup{FSI}$, then $\sbA\cong\mathbf{2}$ or $\sbA$ is anti-idempotent.  (This follows easily from a characterization of FSI De Morgan monoids in \cite[Remark~5.19]{MRW19}, using the fact that $\sbS_3$ cannot be a subalgebra of a homomorphic image of $\sbA$.)  In both cases, $\sbA$ satisfies $f\wedge\neg f^2\leqslant x$, so the statement of Lemma~\ref{dmm lem 1} can be sharpened as follows: in any member of a Koll\'{a}r variety of De Morgan monoids, $f\wedge\neg f^2$ is the least element, whence $e\vee f^2$ is the greatest element.

This fact shows, with the help of (\ref{ddt for rt}), that in any axiomatic extension of $\mathbf{R^t}$ that possesses an IL, the IL can be assumed to take the following form:
\[
\textup{$\Gamma\cup\{\al_1,\dots,\al_n\}$ \,is inconsistent iff\, $\Gamma\,\vdash (\al_1\wedge{\dots}\wedge\al_n\wedge e)\to (f\wedge\neg f^2)$.}
\]
In other words, an IL-sequence $\{\Psi_n:n\in\mathbb{N}^+\}$ for the extension is given by
\[
\Psi_n=\{(v_1\wedge{\dots}\wedge v_n\wedge e)\to (f\wedge\neg f^2)\}.
\]


\begin{thebibliography}{10}

\bibitem{AB75}
    A.R.\ Anderson, N.D.\ Belnap, Jnr., `Entailment: The Logic of Relevance
    and Necessity, Vol.~1', Princeton
    University Press, 1975.



\bibitem{BH06}
    W.J.\ Blok, E.\ Hoogland, {\em The Beth property in algebraic logic},
    Studia Logica {\bf 83} (2006), 49--90.


\bibitem{BP86}
    W.J.\ Blok, D.\ Pigozzi,
    {\em Protoalgebraic logics}, Studia Logica
    {\bf 45} (1986), 337--369.

\bibitem{BP88}
    W.J.\ Blok, D.\ Pigozzi,
    {\em Local deduction theorems in algebraic
    logic},
    in
    H.\ Andr\'{e}ka, J.D.\ Monk, I.\ Nemeti (eds.),
    `Algebraic Logic',
    Colloquia Mathematica Societatis J\'{a}nos Bolyai 54,
    Budapest (Hungary), 1988,
    pp.\,75--109.

\bibitem{BP89}
    W.J.\ Blok, D.\ Pigozzi,
    `Algebraizable Logics', Memoirs of the
    American Mathematical Society~396,
    Amer.\ Math.\ Soc.,
    Providence, 1989.


\bibitem{BP}
    W.J.\ Blok, D.\ Pigozzi,
    {\em Abstract algebraic logic and the
    deduction theorem}, manuscript, 1997.
    [See \verb+http://orion.math.iastate.edu/dpigozzi/+ for updated version, 2001.]



\bibitem{CR17}
    M.A.\ Campercholi, J.G.\ Raftery, {\em Relative congruence formulas and decompositions in quasivarieties},
    Algebra Universalis {\bf 78}
    (2017), 407--425.

\bibitem{CV12}
    M.A.\ Campercholi, D.J.\ Vaggione, {\em Implicit definition of the quaternary discriminator},
    Algebra Universalis {\bf 68}
    (2012), 1--16.

\bibitem{CZ92}
A.\ Chagrov, M.\ Zakharyashchev, \emph{Modal companions of intermediate propositional logics},
Studia Logica \textbf{51} (1992), 49--82.

\bibitem{CN13}
P.\ Cintula, C.\ Noguera,
\emph{The proof by cases property
and its variants in structural
consequence relations},
Studia Logica \textbf{101} (2013), 713--747.

\bibitem{Cze85}
    J.\ Czelakowski,
    {\em Algebraic aspects of deduction theorems},
    Studia Logica {\bf 44} (1985), 369--387.

\bibitem{Cze86}
    J.\ Czelakowski,
    {\em Local deduction theorems},
    Studia Logica {\bf 45}
    (1986), 377--391.

\bibitem{Cze01}
   J.\ Czelakowski,
   `Protoalgebraic Logics', Kluwer, Dordrecht, 2001.


\bibitem{CP99}
    J.\ Czelakowski, D.\ Pigozzi,
    {\em Amalgamation and interpolation in abstract algebraic logic}, in
    X.\ Caicedo, C.H.\ Montenegro (eds.),
    `Models, Algebras and Proofs',
    Lecture Notes in Pure and Applied Mathematics, No.~203,
    Marcel Dekker, New York, 1999, pp.\,187--265.

\bibitem{DL59}
M.A.E.\ Dummett, E.J.\ Lemmon, \emph{Modal logics between $S4$ and $S5$}, Zeitschrift f\"{u}r mathematische Logik und Grundlagen der Mathematik \textbf{5} (1959), 250--264.

\bibitem{Dun66}
    J.M.\ Dunn, `The Algebra of Intensional Logics', PhD thesis, University of Pittsburgh, 1966.

\bibitem{Fon16}
    J.M.\ Font, `Abstract Algebraic Logic -- An Introductory Textbook', Studies in Logic 60,
    College Publications, London, 2016.


\bibitem{Fri62}
O.\ Frink, \emph{Pseudo-complements in semilattices}, Duke Math.\ J.\ \textbf{29} (1962), 505--514.


\bibitem{Gab81}
D.M.\ Gabbay, `Semantical Investigations in Heyting's Intuitionistic Logic', Vol.\,148 of
Synthese Library, D. Reidel, Dordrecht, Boston, 1981.

\bibitem{Gor86}
    V.A.\ Gorbunov, \emph{The cardinality of subdirectly irreducible systems in quasivarieties},
    Algebra and Logic {\bf 25}
    (1986), 1--34.

\bibitem{Gor98}
    V.A.\ Gorbunov,
    `Algebraic Theory of Quasivarieties', Consultants
    Bureau, New York, 1998.





\bibitem{Haj98}
    P.\ H\'{a}jek, `Metamathematics of Fuzzy Logic', Kluwer, Dordrecht, 1998.


\bibitem{Jan68}
V.A.\ Jankov, \emph{Calculus of the weak law of the excluded middle},
Rossiiskaya Akademiya Nauk.\ Izvestiya Seriya Matematicheskaya \textbf{32} (1968), 1044--1051 (Russian).





\bibitem{Kol79}
    J.\ Koll\'{a}r, {\em Congruences and one-element subalgebras}, Algebra Universalis {\bf 9} (1979), 266--267.

\bibitem{KK06}
T.\ Kowalski, M.\ Kracht, \emph{Semisimple varieties of modal algebras}, Studia Logica \textbf{83} (2006), 351--363.

\bibitem{Kra07}
M.\ Kracht, \emph{Modal consequence relations}, in:
P.\ Blackburn, J.\ van Benthem, F.\ Wolter (eds.),
`Handbook of modal logic.'
Studies in Logic and Practical Reasoning, Vol.\,3,
Elsevier B.V., Amsterdam, 2007, pp.\,491–545.

\bibitem{LP}
T.\ L\'{a}vi\v{c}ka, A.\ P\v{r}enosil, \emph{Semisimplicity, the excluded middle and Glivenko Theorems}, manuscript.


\bibitem{MR74}
L.\ Maksimova, V.V.\ Rybakov, \emph{On the lattice of normal modal logics},
Algebra and Logic \textbf{13} (1974), 105--122 (1975).

\bibitem{Mey86}
    R.K.\ Meyer, {\em Sentential constants in R and R$^\neg$}, Studia Logica {\bf 45}
    (1986), 301--327.

\bibitem{MDL74}
    R.K.\ Meyer, J.M.\ Dunn, H.\ Leblanc, {\em Completeness of relevant
    quantification theories}, Notre Dame J.\ Formal
    Logic
    {\bf 15}
    (1974), 97--121.

\bibitem{MRW19}
    T.\ Moraschini, J.G.\ Raftery, J.J.\ Wannenburg, {\em Varieties of De Morgan monoids:
    minimality
    and irreducible algebras}, J.\ Pure Appl.\ Algebra
    {\bf 223}
    (2019), 2780--2803.



\bibitem{MRW20b}
    T.\ Moraschini, J.G.\ Raftery, J.J.\ Wannenburg, {\em Varieties of De Morgan monoids: covers of
    atoms}, Rev.\ Symbolic Logic \textbf{13}
    (2020), 338--374.

\bibitem{MRW20}
    T.\ Moraschini, J.G.\ Raftery, J.J.\ Wannenburg, \emph{Singly generated quasivarieties and residuated structures},
    Math.\ Logic Quarterly \textbf{66}
    (2020), 150--172.

\bibitem{MRW20a}
    T.\ Moraschini, J.G.\ Raftery, J.J.\ Wannenburg, {\em Epimorphisms, definability and cardinalities},
    Studia Logica {\bf 108} (2020),
    255--275.

\bibitem{OR07}
    J.S.\ Olson, J.G.\ Raftery, {\em Positive Sugihara monoids},
    Algebra Universalis {\bf 57} (2007), 75--99.

\bibitem{Pig88}
    D.\ Pigozzi, {\em Finite basis theorems for relatively congruence-distributive
    quasivarieties}, Trans.\ Amer.\ Math.\ Soc.\ {\bf 310}
    (1988), 499--533.

\bibitem{Raf11}
     J.G.\ Raftery, {\em Contextual deduction theorems}, Studia Logica {\bf 99}
     (2011), 279--319.

\bibitem{Raf13}
    J.G.\ Raftery, {\em Inconsistency lemmas in algebraic logic}, Math.\ Logic Quarterly {\bf 59}
    (2013), 393--406.


\bibitem{Ryb97}
    V.V.\ Rybakov, `Admissibility of Logical Inference Rules',
    Studies in Logic and the Foundations of
    Mathematics~136, Elsevier,
    Amsterdam, 1997.



\bibitem{Woj88}
    R.\ W\'{o}jcicki,
    `Theory of Logical Calculi',  Kluwer,
    Dordrecht, 1988.

\end{thebibliography}
\end{document}